\documentclass[11pt]{article}

\usepackage{latexsym,amsfonts,amsmath,amssymb}
\usepackage[english]{babel}
\usepackage[latin1]{inputenc}
\usepackage[T1]{fontenc}
\usepackage{graphics}

\begin{document}

\renewcommand{\baselinestretch}{1.3}
\newtheorem{duge}{Lemma}[section]
\newtheorem{rem}[duge]{Remark}
\newtheorem{prop}[duge]{Proposition}
\newtheorem{defi}[duge]{Definition}
\newtheorem{theo}[duge]{Theorem}
\newtheorem{nota}[duge]{Notation}
\newtheorem{rappel}[duge]{Rappel}
\newtheorem{hypo}[duge]{Hypothesis}
\newtheorem{obj}[duge]{Objectif}

\newcommand{\mcp}{\mathbb{P}}
\newcommand{\mcq}{\mathbb{Q}}
\newcommand{\mce}{\mathbb{E}}
\newcommand{\taux}{\lambda_{b,k}}
\newcommand{\pn}{\mathcal{P}_{n}}
\newcommand{\mcn}{\mathbb{N}}
\newcommand{\pinf}{\mathcal{P}_{\infty}}
\newcommand{\partn}{\{ 1,\dots,n \}}
\newcommand{\sr}{\mathcal{S}}
\newcommand{\srb}{\overline{\mathcal{S}}}
\newcommand{\pic}{\widehat{\Pi}}
\newcommand{\un}{\mathbf{1}}
\newcommand{\limitr}{\underset{s \rightarrow 0}{\lim}}
\newcommand{\N}{\mathbb{N}}
\newcommand{\R}{\mathbb{R}}
\newcommand{\Rp}{\mathbb{R_{+}}}
\newcommand{\pa}{\mathbb{P}_{\alpha}}
\newcommand{\po}{\left}
\newcommand{\pf}{\right}

\title{Ruelle's probability cascades seen as a fragmentation process
}

\author{
Anne-Laure Basdevant}
\date{}
\maketitle

\begin{center}
\it{Laboratoire de Probabilités et Modèles Aléatoires,\\
 Universit\'e Pierre et Marie Curie,\\ 175 rue du Chevaleret,
75013
Paris, France.}
\end{center}

\vspace*{0.8cm}

\begin{abstract}
 In this paper, we study Ruelle's
probability
cascades \cite{Ruelle87} in
the framework of time-inhomogeneous fragmentation processes.
We describe
 Ruelle's cascades mechanism exhibiting a family of
measures  $(\nu_t,t\in [0,1[)$ that characterizes its infinitesimal evolution. To this end, we
will first extend the  time-homogeneous fragmentation theory to
the inhomogeneous case. In the
last section, we will study the behavior for small and large times
of Ruelle's fragmentation process.
\end{abstract}


\bigskip
\noindent{\bf A.M.S. Classification. } 60 J 25, 60 G 09.

\bigskip
\noindent{\bf e-mail. } Anne-Laure.Basdevant@ens.fr

\vspace*{0.8cm}

\section{Introduction}
Ruelle \cite{Ruelle87} introduced a cascade of random probability
measures in  order to study Derrida's GREM model in statistical
mechanics. This approach was further developed by  Bolthausen and
Sznitman \cite{Bolthausensznitman98}, who pointed out that an
exponential time-reversal transforms Ruelle's probability cascades
into a remarkable coalescent process. Previously Neveu
(unpublished) observed that Ruelle's probability cascades were
also related to the genealogy of some continuous state branching
process ; we refer to \cite{Bertoinlegall00} for precise
statements and the connexion with Bolthausen-Sznitman coalescent.
Furthermore, Pitman  \cite{Pitman02} obtained a number of explicit
formulas on the law of Ruelle's cascades ; in particular he showed
that the latter can be viewed as a fragmentation process and
specified its semi-group in terms of certain Poisson-Dirichlet
distributions. Returning to applications to Derrida's GREM model,
we mention the important works by Bovier and Kurkova
\cite{BovierKurkova03,BovierKurkova041,BovierKurkova044} who
established in particular properties of the limiting Gibbs
measure.

The purpose of this paper is to dwell on Pitman's observation that
Ruelle's cascade can be viewed as a time-inhomogeneous
fragmentation process. The  theory of time-homogeneous
fragmentation processes was developed recently (see eg
\cite{Berestycki02,CoursBertoin03,Bertoin01}), and we shall
briefly show how it can be extended to the time-inhomogeneous
setting. Roughly the basic result is that the distribution of a
time-inhomogeneous fragmentation can be characterized by a
so-called instantaneous rate of erosion
 (which is a non-negative real number that depends on the time parameter),
 and an instantaneous dislocation measure (which specifies the rate of sudden dislocation).
We shall establish that for Ruelle's probability cascades, the
instantaneous erosion is zero, and we will provide several
descriptions of the instantaneous dislocation measure.
Specifically, the latter is related to the well-known
Poisson-Dirichlet distributions, in particular we shall establish
a stick-breaking construction, compute the corresponding
exchangeable partition probability function, and derive some
relations of absolute continuity. In this direction, we mention
that related (but somewhat less precise) results have been proven
independently by Marchal \cite{Marchal2}. Finally, as examples of
applications, we shall prove some asymptotic results for Ruelle's
probability cascades at small and large times.

The rest of this work is organized as follows. The next section is
devoted to preliminaries, then we briefly present the extension of
the theory of fragmentation processes to the time-inhomogeneous
setting. The main results on Ruelle's probability cascades are
established in section \ref{appl}, and finally section deals with
applications to the asymptotic behavior.

\section{Preliminaries}
\subsection{Ruelle's cascades and their representation with stable subordinators}\label{cascades}

Let us briefly recall  the construction of Ruelle's cascades \cite{Bertoinlegall00,Bolthausensznitman98,Ruelle87}.
Let $p>1$  be an integer and  let $0<x_1<\ldots<x_p<1$ be a finite sequence of real numbers.
For $k\in \{1,\dots,p\}$,  $(\eta_{i_1,\ldots,i_k},i_1 \dots i_k\in\mcn)$ denotes a family of random variables such that :
\begin{itemize}
\item for $k\in \{1,\dots,p\}$, $i_1,\ldots,i_{k-1}\ge 1$ fixed,  the distribution of $(\eta_{i_1,\ldots,i_{k-1},j},j\in\mcn)$ is that of the sequence of atoms of
a Poisson measure on $]0,\infty[$ with intensity  $x_kr^{-1-x_k}dr$, arranged according to the decreasing order of their sizes,
\item the families $(\eta_{i_1,\ldots,i_{k-1},j},j\in\mcn)$ for $k\in \{1,\dots,p\}$, $i_1,\ldots,i_{k-1}\ge 1$ are independent.
\end{itemize}

Set $\theta_{i_1,\ldots,i_k}=\eta_{i_1}\ldots\eta_{i_1,\ldots,i_k}$. We can easily show that $C=\sum_{i_1\dots i_p}\theta_{i_1,\ldots,i_p}$ is
 almost surely finite.
Next we define Ruelle's cascades :
$$\overline{\theta}_{i_1,\ldots,i_p}=\frac{\theta_{i_1,\ldots,i_p}}{C}\hspace*{1cm} \mbox{ and recursively }\hspace*{1cm} \overline{\theta}_{i_1,\ldots,i_{k-1}}=
\sum_{j=1}^{\infty}\overline{\theta}_{i_1,\ldots,i_{k-1},j}. $$

Bertoin and Le Gall \cite{Bertoinlegall00} have proved we can relate this  process to the genealogy of  Neveu's CSBP (continuous-state branching process).
Precisely, they have proved that there exists  a process $(S^{(s,t)}(a),0\le s<t,a\ge 0)$ such that :
\begin{itemize}
\item $\forall 0\le s<t,\; S^{(s,t)}= (S^{(s,t)}(a),a\ge 0)$ is a
stable subordinator with index $e^{-(t-s)},$ \item $\forall p\ge
2,\; 0\le t_1 \le \ldots \le t_p,\; S^{(t_1,t_2)},\ldots,
S^{(t_{p-1},t_p)}$ are independent and
$S^{(t_1,t_p)}(a)=S^{(t_{p-1},t_p)}\circ\ldots\circ
S^{(t_1,t_2)}(a).$
\end{itemize}

Set  $0<t_1<\ldots<t_p$ such that

$$x_1=e^{-t_p} \mbox{ and } x_k=e^{-(t_p-t_{k-1})}, k=2,\ldots,p. $$
Let us fix $a>0$. We define recursively, for $k=1,\ldots,p$, random intervals $D^{(t_1,\ldots,t_k,a)}_{i_1,\ldots,i_k}$ in the following way :\\
$D^{(a)}=]0,a[$.\\
Let $k\ge 1$, $i_1,\ldots,i_{k-1}\in \mcn$. Let $(b_{i_1,\ldots,i_{k}},i_k\in\mcn)$ be the jump times of  $S^{(t_{k-1},t_k)}$ on the interval
 $D^{(t_1,\ldots,t_{k-1},a)}_{i_1,\ldots,i_{k-1}}$ listed in the decreasing order of sizes. We set
\begin{equation}\label{inter}
D^{(t_1,\ldots,t_k,a)}_{i_1,\ldots,i_k}=]S^{(t_{k-1},t_k)}(b_{i_1,\ldots,i_{k}}-),S^{(t_{k-1},t_k)}(b_{i_1,\ldots,i_{k}}) [ \; \mbox{ and } \;
\xi^{(t_1,\ldots,t_k,a)}_{i_1,\ldots,i_k}=|D^{(t_1,\ldots,t_k,a)}_{i_1,\ldots,i_k}|.
\end{equation}

Bertoin et Le Gall have proved that the families
$$\left(\left(S^{(0,t_p)}(a)\right)^{-1} \xi_{i_1,\ldots,i_p}; i_1,\ldots,i_p\in \mcn \right) \mbox{ and }
 \left(\overline{\theta}_{i_1,\ldots,i_p}; i_1,\ldots,i_p\in \mcn\right)$$
have the same law.

\subsection{Ruelle's cascades as fragmentation processes}
Using this representation of Ruelle's cascades in terms of stable subordinators, we can exhibit a link with fragmentation processes.

Recall that the law $\beta\left(a,b\right)$ has density
$$\frac{\Gamma\left(a+b\right)}{\Gamma\left(a\right)\Gamma\left(b\right)}x^{a-1}\left(1-x\right)^{b-1}\mathbf{1}_{[0,1]}dx,$$
and  let us introduce some definition :

\begin{defi}\label{PD}\cite{Pitmanyor97}
For $0\le \alpha \le 1$, $\theta>-\alpha$, let
$\left(Y_n\right)_{n\ge 1}$ be a sequence of independent random variables with respective laws
$\beta\left(1-\alpha,\theta+n\alpha\right)$. Set
\begin{equation*}
\widehat f_1=Y_1 \hspace*{0.6cm} \widehat
f_n=\left(1-Y_1\right)\dots\left(1-Y_{n-1}\right)Y_n
\hspace*{0.6cm} \widehat f=\left(\widehat f_n\right)_{n\ge 1}.
\end{equation*}
Then $\sum_{i}{\widehat f_{i}}= 1$. Let
$f=\left(f_n\right)_{n>0}$ be the decreasing rearrangement of the sequence
 $(\widehat f_{n})_{n\ge 1}$.We define the Poisson-Dirichlet law with parameter
$\left(\alpha,\theta\right)$, denoted
 $PD\left(\alpha,\theta\right)$,
 as the distribution of  $f$.
\end{defi}

Thereafter $\sr$ stands for the set of decreasing sequences of
non-negative numbers with sum  equal to $1$.
 $\sr$  is endowed with the uniform distance. $\srb$ denotes its closure, it is the set
 of decreasing sequences of non-negative real numbers whose  sum is less than or equal to $1$ and is called the set of mass-partitions.
Notice that $\srb$ is a compact set.

\begin{defi} Let $s=(s_i,i\in \mcn)$ be an element of $\srb$ and
$s^{(.)}=(s^{(i)},i \in \mcn)$  a sequence in $\srb$. Consider the fragmentation of $s_i$ by
$s^{(i)}$, i.e. the sequence $\tilde{s}^{(i)}=(s_is^{(i)}_j,j\in \mcn)$.
 The  decreasing rearrangement of
all the terms
 of the sequences  $\tilde{s}^{(i)}$ as $i$ describes $\mcn$ is called fragmentation of $s$ by
 $s^{(.)}$.
If $\mathbb{P}$ is a probability on $\srb$, we define the transition kernel
 $\mathbb{P}-FRAG\left(s,.\right)$ as the distribution of a fragmentation of
 $s$ by $s^{(.)}$, where $s^{(.)}$ is an iid sequence of random mass-partition with law $\mathbb{P}$.
\end{defi}

A Markov process  $(F(t),t\in [0,1[)$ with values in $\srb$
is called a fragmentation process if the following properties are fulfilled :
\begin{itemize}
\item $F(t)$ is continuous in probability.
\item Its semi-group has the following form :\\
 for all $t,t'\in  [0,1[$ such that $t+t'\in[0,1[$, the conditional
law of $F(t+t')$ given $F(t)=s$ is the law of
$\mathbb{P}_{t,t+t'}-FRAG(s,\cdot)$ where $\mathbb{P}_{t,t+t'}$ is a probability on
$\sr$.
\end{itemize}
 The fragmentation is said homogeneous (in time) if
$\mathbb{P}_{t,t+t'}$
 depends only on $t'$.
Besides,  $(F(t),t\in [0,1[)$ is called a standard fragmentation
process if $F(0)$ is almost surely equal to the sequence
$\mathbf{1}=(1,0,\ldots)$.

In the case of Ruelle's cascades, using the work of Bertoin et Pitman \cite{Bertoin00} (Lemma 9), we know that for any integer  $2\le k\le p$,
$(\overline{\theta}_{i_1,\ldots,i_k},i_1,\ldots,i_{k}\ge 1)$ is a $PD(x_{k},-x_{k-1})$-fragmentation of
$(\overline{\theta}_{i_1,\ldots,i_{k-1}},i_1,\ldots,i_{k-1}\ge 1)$.
More precisely we have :

\begin{prop}\label{frag}
 There exists  a   time-inhomogeneous fragmentation
 $(F(t),t\in[0,1])$ with semi-group  $\mathbb{P}_{t,t+t'}=PD(t+t',-t)$ such that
$$\Big(\left(\overline{\theta}_{i_1}; i_1\in \mcn\right),\ldots \left(\overline{\theta}_{i_1,\ldots,i_p}; i_1,\ldots,i_p\in \mcn\right)\Big)
\overset{law}{=} \Big(F(x_1),\ldots,F(x_p)\Big).$$
\end{prop}

 In the sequel, we  call $F$, Ruelle's fragmentation.
To study Ruelle's cascade, it should be possible to use the
fragmentation process theory developed for example in
\cite{Bertoin03}, but first, we must extend this theory to
time-inhomogeneous fragmentations.

\subsection{Exchangeable random partitions}
In this section, we recall the connections between exchangeable
random partitions and
mass-partitions. Let us first introduce some useful notation :\\
we denote by $\mcn$  the set of positive integers. For $n\in \mcn$, $[n]$ denotes the set $\{1,\ldots,n\}$ and
$\pn$ denotes the set of partitions of $[n]$, $\pinf$
the set of partitions of $\mcn$. For all $n<m$, for all $\pi \in \mathcal{P}_m,\;
\pi_{|n}$ denotes the restriction of $\pi$ to $\pn$.
We endow $\pinf$  with the distance
$d(\pi,\pi')=\frac{1}{\sup\{n\in \mcn\; \pi_{|n}=\pi'_{|n}\}}$.
 The partition with a single block is denoted by $\un$.
 We always label the blocks of a partition according to the increasing
order of their smallest element.

A random partition of  $\mcn$ is called exchangeable if its distribution is invariant by the action of the group of finite permutations of $\mcn$.
Kingman \cite{Kingman82} has proved that each block of an exchangeable random partition has a frequency, i.e.
if
 $\pi=(\pi_1,\pi_2,\ldots)$ is  an exchangeable random partition, then
 $$\forall i \in \mcn \hspace*{1cm}
f_i=\underset{n \rightarrow \infty}{\lim}\frac{\sharp\{\pi_i\cap
[n]\}}{n} \hspace*{1cm}\mbox{exists a.s.}$$
One calls
 $f_i$  the frequency of the block $\pi_i$.
Therefore, for all exchangeable random partitions, we can associate a probability on $\srb$ which will be the law of the decreasing rearrangement of the sequence of
the partition frequencies.

Conversely, given a law $\mathbb{P}$ on $\srb$, we can construct
an exchangeable random partition whose law  of its  frequency  sequence is $\mathbb{P}$ (cf. \cite{Kingman82}).
Let us specify this construction :
we pick $s\in \srb$ with law $\mathbb{P}$ and we draw a sequence of independent random variables $U_i$ with uniform law on $[0,1]$.
Conditionally on $s$, two integers $i$ and $j$ are in the same block of $\Pi$ iff there exists an integer $k$ such that
$\sum_{l=1}^{k}{s_l}\le U_i<  \sum_{l=1}^{k+1}{s_l} $ and
$\sum_{l=1}^{k}{s_l}\le U_j<  \sum_{l=1}^{k+1}{s_l} $.
This construction of a law on the set of partitions from a law on $\srb$ is often called  ``paint-box process''.

Kingman's representation Theorem states that any  exchangeable
random partition can be constructed in this way. Therefore, we
have a natural bijection between the laws on $\srb$ and the laws
on  the exchangeable random partitions. \\
We also define an exchangeable measure  $\rho_\nu$ on $\pinf$ from
a measure $\nu$ on $\srb$ by :
$$\rho_{\nu}(\cdot)=\int_{\srb}{\rho_u(\cdot)\nu(du)}$$
where $\rho_u$ is the law on $\pinf$ obtained by the paint-box
based on the mass-partition $u$.

For any  exchangeable random partition $\Pi$, we define a
symmetric function $p$ on finite sequences of $\mcn$ such that,
 for every  $n,n_1,\ldots,n_k$  integers with  $n=n_1+\ldots+n_k$, $$p(n_1,\ldots,n_k)=\mcp(\Pi_{|n}=\pi),$$ where $\pi$ is a partition of
$[n]$ with $k$ blocks of size $n_1,\ldots,n_k$.
The fact that $\mcp(\Pi_{|n}=\pi)$  depends only on $n_1,\ldots,n_k$ stems from the  exchangeability of $\Pi$. One calls $p$ the EPPF (exchangeable
partition probability function) of $\Pi$.

\begin{prop}\cite{Pitman96,Pitman99}\label{alpha-theta}
Let $\widehat f=\left(\widehat f_n\right)_{n\in\mcn}$ be a sequence of random variables
of $\lbrack 0,1\rbrack$
defined as in  Definition \ref{PD} .
Then there exists an exchangeable random partition with frequency distribution $\widehat f$, where $\widehat f_i$ is the $i$-th block frequency
 and where the blocks are listed order of their smallest element. It is a $\left(\alpha,\theta\right)$-partition.
Besides  the EPPF of this partition  is
\begin{equation}\label{alpheq}
p_{\alpha,\theta}\left(n_1,\dots,n_k\right)=\frac{\lbrack
\frac{\theta}{\alpha}\rbrack_k}{\lbrack
\theta\rbrack_n}\prod_{i=1}^{k}{-\lbrack -\alpha\rbrack_{n_i}} \hspace*{0.5cm}\mbox{for } \theta\neq 0
 \hspace*{2cm}\mbox{(Ewens-Pitman's formula)}
\end{equation}
where $\lbrack x\rbrack_n=\prod_{i=1}^{n}\left(x+i-1\right)$ and $n=\sum_{i=1}^{k} {n_i}$.\\
For $\theta=0$, the formula is extended by continuity.
\end{prop}

This proposition also proves that the law of the sequence  $\widehat
f$ is invariant by  size-biaised rearrangement.

\vspace*{0.5cm}

In the case of Ruelle's fragmentation, we know that,  at time $t$,
$F(t)$ has the $PD(t,0)$ law. So we have the following proposition :

\begin{prop}\label{EPFBS}The EPPF $q_t$ of the random partition  associated with Ruelle's fragmentation at time $t$, $F(t)$, is :
\begin{equation}\label{EPFBSF}
q_{
t}\left(n_1,\ldots,n_k\right)=\frac{\left(k-1\right)!}
{\left(n-1\right)!}t^{k-1}\prod_{i=1}^{k}{[1-t]_{n_i-1}}
\end{equation}
where $n=\sum_{i=1}^{k} n_i.$
\end{prop}

\begin{rem}
We can also construct a random partition with distribution $p_{\alpha,\theta}$ recursively
(Chinese restaurant construction):\\
First, the integer $1$ necessarily belongs to the first block, denoted $B_1$.
Suppose  the $n$ first integers split up in $b$ blocks : $\Pi_n=\left(B_1,\ldots,B_b\right)$,
where  block $B_i$ has cardinal $n_i$. We now define
$\Pi_{n+1}$ with the following rule :\\
$\mcp\left(\Pi_{n+1}=(B_1,\ldots,B_i\cup\{n+1\},\ldots,B_b)\right)=\frac{n_i-\alpha}{n+\theta}$\\
$\mcp\left(\Pi_{n+1}=(B_1,\ldots,B_b,\{n+1\})\right)=\frac{b\alpha+\theta}{n+\theta}$. \\
Then $\Pi$ is a  $(\alpha,\theta)$-partition
(cf. \cite{Pitman95}).
\end{rem}

\vspace*{0.5cm}

We can also define a notion of fragmentation process  of exchangeable partitions such that there is still a bijection with
 fragmentation processes  of mass-partitions :

 Set $A\subseteq B \subseteq \mcn$ and $\pi \in
\mathcal{P}_{A}$ with $\# \pi=n$. Let $\pi^{(.)}=(\pi^{(i)},i \in
\partn)$, $\pi^{(i)} \in \mathcal{P}_{B}$ for all $i$. Consider
the partition of the $i$-th block of $\pi$,$\pi_i$,
induced by
$\pi^{(i)}$, i.e. $\pi^{(i)}_{|\pi_i}=\tilde{\pi}^{(i)}$.\\
As $i$ describes  $\partn$, the blocks of $\tilde{\pi}^{(i)}$
form the blocks of a partition $\tilde{\pi}$ of $A$. This partition is
 denoted $FRAG(\pi,\pi^{(.)})$. This is the  fragmentation of $\pi$
by $\pi^{(.)}$. \\
If $\mathbb{P}$ is a probability on $\mathcal{P}_{B}$, define the
transition kernel $\mathbb{P}-FRAG\left(\pi,.\right)$ as the distribution of a
fragmentation of $\pi$ by $\pi^{(.)}$, where $\pi^{(.)}$ is a sequence of
iid partition with law $\mathbb{P}$.\\
Let $(\Pi(t),t \in  [0,1[)$ be a  Markov process on
$\pinf$. We call $(\Pi(t),t \in  [0,1[)$  an exchangeable fragmentation process if
the following properties are fulfilled :
\begin{itemize}
\item $\Pi(t)$ is continuous in probability.
\item Its semi-group has the following form  :\\
for all $ t,t^{'}\ge 0$ such that $t+t'<1$, the conditional law of $\Pi(t+t^{'})$
given
$ \Pi(t)=\pi$ is  $\mathbb{P}_{t,t'}-FRAG(\pi,\cdot)$, where $\mathbb{P}_{t,t'}$ is an exchangeable probability on $\pinf$.
\end{itemize}
The fragmentation is homogeneous if $\mathbb{P}_{t,t'}$ depends
only on $t'$.
Furthermore, $(\Pi(t),t\in [0,1[)$ is a standard fragmentation process if $\Pi(0)$ is equal to  $\mathbf{1}$.

\vspace*{0.4cm}

We can check that, with these definitions, if
 $(\Pi(t),t \in [0,1[
)$ is a fragmentation
 process on
partitions, then $(F(t),t\in  [0,1[)$ the frequency process of
$\Pi$, is a  fragmentation
 process on
mass-partitions. Furthermore, the converse is true, i.e., if one
takes a fragmentation process on mass-partitions, then one can
construct a fragmentation
 process on
partitions $\Pi$ such that the frequency process of $\Pi$ is equal
to the initial fragmentation process (cf. \cite{Berestycki02}).

We also remark  that if we consider a fragmentation $(\Pi(t),t\in
[0,1[)$ with semi-group $\mathbb{P}_{t,t'}-FRAG$, then its
restriction to $\pn$, $(\Pi_{|n}(t),t\in [0,1[)$, is a Markov
process with semi-group $\mathbb{P}_{t,t'}^n-FRAG$ where
$\mathbb{P}_{t,t'}^n$ is equal to $\mathbb{P}_{t,t'}$ restricted
to $\pn$ (cf. \cite{CoursBertoin03}).

\vspace*{0.5cm}

For Ruelle's fragmentation, we have an explicit construction of
its corresponding fragmentation on the partitions. Indeed, recall
the representation of Ruelle's cascades with the jumps of a family
of subordinators (cf. Section \ref{cascades}). Let
$(\sigma^{*}_t,t\in[0,1[)$ be a family of stable subordinators
such that for every $0\le t_p<\ldots<t_1<1$, the joint
distribution of $\sigma^{*}_{t_1},\ldots,\sigma^{*}_{t_p}$ is the
same as that of  $\sigma_{t_1},\ldots,\sigma_{t_p}$ with
$\sigma_{t_i}=\tau_{\beta_1}\ldots\tau_{\beta_i}$ where
$t_i=\beta_1\ldots\beta_i$
 and $\tau_{\beta_1},\ldots,\tau_{\beta_p}$ are independent stable subordinators with indices $\beta_1,\ldots,\beta_p$.
For $t\in ]0,1[$,
let $M_t$ be the closure of $\{\sigma^{*}_{t}(u),u\ge 0\}$. Consider
then the family of open subsets of $[0,1[$ :
$G(t)=[0,1[\backslash M_t$, for $t\in [0,1[$. Then $(G(t),t\in  [0,1[)$ is a
nested family, i.e. $G(t)\subset G(s)$ for $0<s<t<1$ and
furthermore, if $F(t)$ is the sequence of ranked lengths of
component intervals of $G(t)$, then $(F(t),t\in [0,1[)$ has the law of
Ruelle's fragmentation ; see
 \cite{Bertoin00}.

Set $0\le t_1<\ldots<t_p<1$. Let us now draw $(U_i)_{i\in \mcn}$, uniform and independent random variables  on $]0,1[$.
For $1\le k \le p$, we construct a partition $\Pi(k)$ of $\mcn$ with the rule :
$$i\overset{\Pi(k)}\sim j \Leftrightarrow U_i \mbox{ and } U_j \mbox{ are in the same component interval of } G(t_i). $$

Then $(\Pi(1),\ldots, \Pi(k))$ has the law of a Ruelle's fragmentation on the partition at times $(t_1,\ldots, t_p)$.

\vspace*{0.5cm}

\subsection{Connexion with  Bolthausen-Sznitman's coalescent}
Bolthausen et Sznitman \cite{Bolthausensznitman98} have shown that
it is possible to formulate Ruelle's fragmentation as a coalescent process if we reverse time.
Moreover, for a good choice for the time reversal, the coalescent process is time-homogeneous \cite{Bolthausensznitman98}.
Let us first recall the definition of a coalescent process.

\vspace*{0.2cm}

Set $s\in\srb$ and  let $\Pi=\{B_1,B_2,\ldots\}$ be a partition of
$\mcn$. Set $\tilde{s}_i=\sum_{j\in B_i}{s_j}$.
 The $\Pi$-coagulation of
$s$, denoted $COAG(s,\Pi)$ is the decreasing rearrangement of the sequence $(\tilde{s}_i,i \in \mcn)$. \\
If $\mathbb{P}$ is a probability on  $\sr$, we define the transition kernel $\mathbb{P}-COAG\left(s,.\right)$ as the distribution of a
$\Pi$-coagulation of $s$, where $\Pi$ has the law on $\pinf$ obtained from $\mathbb{P}$ by the paint-box construction.\\
Let $(C(t),t\ge 0)$ be a Markov process on $\srb$.
$(C(t),t\ge 0)$ is a time-homogeneous mass-coalescent process if the following properties are fulfilled :
\begin{itemize}
\item $C(t)$ is continuous in probability.
\item Its semi-group has the following form :\\
for all $t,t'\ge 0$, the conditional law of  $C(t+t')$ given
$C(t)=s$ is the law of  $\mathbb{P}_{t'}-COAG(s,\cdot)$ where
$\mathbb{P}_{t'}$ is a probability on $\sr$.
\end{itemize}
\vspace*{0.2cm}

To see that Ruelle's fragmentation reversed in time is a time-homogeneous coalescent process, we use the following property :

\begin{prop}\cite{Pitman99}\label{fragcog} Set $\alpha \in \rbrack
0,1\lbrack$, $\beta \in \lbrack 0,1\lbrack$ and
$\theta>-\alpha\beta$. The following assertions are equivalent :
\begin{itemize}
\item  $s$ has  $PD\left(\alpha,\theta\right)$ distribution and
$s'$ is a
$PD\left(\beta,\theta/\alpha\right)$-coagulation of $s$.
\item $s'$ has
 $PD\left(\alpha\beta,\theta\right)$ distribution and $s$ is a $PD\left(\alpha,-\alpha \beta\right)$-fragmentation of $s'$.
\end{itemize}
\end{prop}

Thus, if we define $C(t)=F(e^{-t})$ where $(F(t),t\in[0,1[)$ is
Ruelle's fragmentation,  then $(C(t),t\ge 0)$ is a homogeneous
coalescent process with semi-group $PD(e^{-t},0)$-COAG. This
process is called the Bolthausen and Sznitman's coalescent.

Just like the case of  the fragmentation processes, we can associate  a coalescent process on exchangeable
partition to any
mass-coalescent process. For the  Bolthausen-Sznitman's coalescent process on
 partitions, we have an explicit construction \cite{Pitman99}.
It is a simple exchangeable coalescent process, i.e, at each
jump-time of the process $\Pi_{n}(t)$, only one block can be
formed. The jump rates of this process can be explicitly written.
 If we start from a partition with $b$ blocks, each $k$-uplet of blocks coagulates with  rate $\taux$ that
depends only on $b$ and  $k$ and which is equal to :
\begin{equation*}
\taux=\frac{\left(k-2\right)!\left(b-k\right)!}{\left(b-1\right)!}=\int_{0}^{1}{x^{k-2}\left(1-x\right)^{b-k}dx.}
\end{equation*}

\begin{rem}
One can be surprised that an homogeneous Markov process becomes an
inhomogeneous Markov process after time-reversal. In
fact, Ruelle's fragmentation can also be seen as a homogeneous
Markov process, but, if one takes this point of view, it is no
longer a fragmentation process since the evolution of a
particle depends on the other particles. Actually, it is known
that if a random variable $x=\left(x_1,x_2,\ldots\right) \in \srb$
has the PD$\left(\alpha,0\right)$ law , then $\underset{n \rightarrow
\infty}{\lim}\frac{\alpha \ln x_n}{\ln
n}=-1$ (cf. \cite{Pitmanyor97}).\\
In particular, in the case of Ruelle's fragmentation, $F(t)$  has law
 PD$\left(t,0\right)$, therefore
$$t=-\underset{n \rightarrow \infty}{\lim}\frac{\ln n}{\ln
x_n\left(t\right)}.$$ Let $\left(p_{t,t+s}\right)_{t,s>0}$ be the
transition probabilities of $F$. Suppose that the process is in
state $x \in \srb$. For all  $t\in [0,1[$, the process $F$ has a
Poisson-Dirichlet law, so $T(x)=-\underset{n \rightarrow
\infty}{\lim}\frac{\ln n}{\ln x_n} $ exists and $T(x)$ determines
the considered time. For  $y \in \srb$. We define
$$q_s\left(x,y\right)=p_{T\left(x\right),T\left(x\right)+s}\left(x,y\right).$$
Then $\left(q_s\right)_{s\in [0,1[}$ is a
homogeneous transition kernel for $F$.
However, remark that to determine $T(x)$, we must know the other
particles state and the branching property is lost.
\end{rem}

\section{General theory of time-inhomogeneous fragmentation processes}\label{inh}
In this section, we extend the theory of time-homogeneous
fragmentations to time-inhomogeneous fragmentations. For this, we
will first work on fragmentations of partition and next on mass-fragmentations.

\subsection{Measure of an inhomogeneous fragmentation}
Let us first  define precisely  the class of fragmentations we consider
(which includes Ruelle's fragmentation).
We denote $\pn \setminus \{\un\}$ by $\pn^{*}$.

\begin{hypo}\label{hyp}
In the sequel, we always suppose that $\left(\Pi(t),t\in
[0,1[\right)$ is a standard time-inhomogeneous exchangeable
fragmentation for which the following properties are fulfilled :
\begin{itemize}
\item for all $n \in \mcn$, let $\tau_n$ be the time of the first jump
of  $\Pi_{|n}$ and $\lambda_n$ be its law. We have
$$\forall t\in [0,1[, \; \overline{\lambda}_n(t):=\lambda_n\left([t,1[\right)>0$$
and $\lambda_n$ is absolutely  continuous with respect to Lebesgue
measure with  continuous density $g_n(t)$. \item for all $\pi \in
\pn^*,\hspace*{0.3cm}
h_{\pi}^n(t)=\mcp(\Pi_{|n}(t)=\pi\;|\;\tau_n=t)$ is a continuous
function of t.
\end{itemize}
\end{hypo}

Let us  now define an instantaneous jump rate for a fragmentation fulfilling  Hypothesis \ref{hyp}.\\
Set $\pi \in \pn^*$. Set
$h_{\pi}^n\left(t\right)=\mcp\left(\Pi_{|n}\left(t\right)=\pi \;|\;
\tau_n=t\right)$.
 It is the law of the jump given $\tau_n$.\\
 We set $$f_n\left(t\right)=\limitr\frac{1}{s}\mcp\left(\tau_n \in [t,t+s] \;|
  \; \tau_n\ge t\right)=\frac{g_n\left(t\right)}{\overline{\lambda}_n\left(t\right)},$$

and
$$q_{\pi,t}=h_{\pi}^n\left(t\right)f_n\left(t\right)=\limitr\frac{1}{s}\mcp\left(\Pi_{|n}\left(\tau_n\right)=\pi
\; \& \;\tau_n \in [t,t+s]\;|\; \tau_n\ge t\right).$$

It is the probability density that the process $\Pi_{|n}$ jumps at
time $t$ from the state $\mathbf{1}$ to the state $\pi$ given that
$\Pi_{|n}$ has not jump before.

\begin{prop}\label{taux}
For $\pi \in \pn$, $n'\ge n$, set $Q_{n',\pi}=\{\pi^{'} \in \mathcal{P}_{n'},
\pi^{'}_{|n}=\pi\}$.
For each  $t\in [0,1[$, there exists a unique
measure $\mu_t$ on $\pinf$ such that
$$\forall n \in \mcn \; \forall \pi \in \pn^* \;
\; \mu_t\left(Q_{\infty,\pi}\right)=q_{\pi,t} \mbox{ and }
\mu_t\left(\mathbf{1}\right)=0.$$ The family of measures
$\left(\mu_t,t\ge 0\right)$ characterizes the law of the
fragmentation.
\end{prop}

Proof. We have
\begin{equation}\label{somme}
\forall n'>n, \forall \pi \in \pn^*
 \sum_{\pi' \in Q_{n',\pi}}{q_{\pi',t}}=
q_{\pi,t}.
\end{equation}

In fact, at time $t$, if the process $\Pi_{|n'}$ has not jumped
yet, it will jump between time
 $t$ and time $t+dt$ to the state such that
$\Pi_{|n}=\pi$ with probability $\sum_{\pi' \in
Q_{n',\pi}}{q_{\pi',t}}dt$.
Besides, we have the following equality :
$$\mcp\left(\Pi_{|n}\left(\tau_n\right)=\pi \; \&
\;\tau_n\in [t, t+dt] \;|\; \tau_n\ge
t\right)=\mcp\left(\Pi_{|n}\left(\tau_n\right)=\pi \; \&
\;\tau_n\in [t, t+dt] \;|\; \tau_{n'}\ge t\right),$$ because the
event that the block $[n']$ has already split, does not affect the
process $\Pi_{|n}$. In fact, as $(\Pi_{|n}(t),t\in [0,1[)$ is a
Markov process, the law of the process $(\Pi_{|n}(t),t\in
[t_0,1[)$ depends only on $\Pi_{|n}(t_0)$. Therefore, we have
equality (\ref{somme}).

Let us  now define  $\mu_t(Q_{\infty,\pi})=q_{\pi,t}$. By (\ref{somme}), this application can be extended to an additive application.
 By  Caratheodory's Theorem,  $\mu$ can be extended to an unique measure on $\pinf$.
\vspace*{0.4cm}

We have now to prove that this family of measures determines the
fragmentation law. To this end we just have to prove that the
family of measures $\left(\mu_t, t\in [0,1[\right)$
characterizes every jump rate of $\Pi_{|n}\left(t\right)$.
Set  $\pi,\pi' \in \pn$, $t_0 \in [0,1[$. Let $\tau_n'$ be the time of the first jump
of $\Pi_{|n}\left(t\right)$ after
$t_0$. We must express
$\limitr\frac{1}{s}\mcp\left(\Pi_{|n}\left(\tau_n'\right)=\pi' \; \& \;\tau_n'\in
[t, t+s]
\;|\; \Pi_{|n}\left(t_0\right)=\pi\right)$  in terms of $\left(\mu_t, t\in [0,1[\right)$. \\
If $\pi'$ can not be obtain from the fragmentation of one
block of $\pi$, we clearly have :
$$\mcp\left(\Pi_{|n}\left(\tau_n'\right)=\pi' \; \& \;\tau_n'\in [t, t+dt]
\;|\; \Pi_{|n}\left(t_0\right)=\pi\right)=0.$$ Permuting the
indices
 (which does not change the law by exchangeability), we can suppose
 $\pi=\left(A_1,\ldots,A_N\right)$ and
$\pi'=\left(B_1,\ldots,B_k,A_2,\ldots,A_N\right)$ where $
\pi''=\left(B_1,\ldots,B_k\right) \in \mathcal{P}_{A_1}$. Let
 $\tau_{[A_i]}'$ be the first jump time of $\Pi_{|A_i}$.
Then, by the branching property :
$$\mcp\left(\Pi_{|n}\left(\tau_n'\right)=
 \pi' \;
\& \;\tau_n'\in [t, t+dt] \;|\;
\Pi_{|n}\left(t_0\right)=\pi\right)\hspace*{7cm}$$
 \begin{eqnarray*}
&= &\mcp\left(\Pi_{|A_1}\left(\tau_{[A_1]}\right)=\pi'' \; \&
\;\tau_{[A_1]}\in [t, t+dt] \;|\;
\Pi_{|A_1}\left(t_0\right)=\un\right)\prod_{i=2}^{N}{\mcp\left(\tau_{[A_i]}>t\;|\;\tau_{[A_i]}>t_0\right)}\\
&= &\mcp\left(\Pi_{|A_1}\left(\tau_{[A_1]}\right)=\pi'' \; \&
\;\tau_{[A_1]}\in [t, t+dt] \;|\;
\tau_{[A_1]}>t_0\right)\prod_{i=2}^{N}{\mcp\left(\tau_{[A_i]}>t\;|\;\tau_{[A_i]}>t_0\right)}\\
&= &\mcp\left(\Pi_{|A_1}\left(\tau_{[A_1]}\right)=\pi'' \; \&
\;\tau_{[A_1]}\in [t, t+dt] \;|\;
\tau_{[A_1]}>t_0\right)\prod_{i=1}^{N}{\frac{\overline{\lambda}_{|A_i|}\left(t\right)}{\overline{\lambda}_{|A_i|}\left(t_0\right)}}.
\end{eqnarray*}

Thus we have :

$$\limitr\frac{1}{s}\mcp\left(\Pi_{|n}\left(\tau_n'\right)=\pi' \; \& \;\tau_n'\in
[t, t+s]
\;|\; \Pi_{|n}\left(t_0\right)=\pi\right)
=\mu_t\left(Q_{\infty,\pi''}\right)
\prod_{i=1}^{N}{\frac{\overline{\lambda}_{|A_i|}\left(t\right)}{\overline{\lambda}_{|A_i|}\left(t_0\right)}}$$

and $\overline{\lambda}_{n}$ is easily expressed as a  function of
$\mu_t$ (cf. below).

\begin{prop}\label{intmut}
The application from $[0,1[$ to the set of measures on $\pinf$
which at $t$ associates  $\mu_t$ constructed according to the
 proposition above, verifies :
 \begin{itemize}
 \item $\mu_t$ is  an exchangeable measure
 such that \\
 $\mu_t\{\un\}=0$ and
 $\forall n \in \mcn\;\mu_t\left(\{\pi \in \pinf, \pi_{|n}\neq \mathbf{1}
\}\right)<\infty,$
 \item $\forall n \in \mcn\; \forall t\in [0,1[ \mbox{ we have }\int_{0}^{t}{\mu_u\left(\{\pi \in \pinf, \pi_{|n}\neq \mathbf{1}
\}\right)du}<\infty.$
\end{itemize}
\end{prop}

Proof. The exchangeability is clear and we have
$$\mu_t\left(\{\pi \in \pinf, \pi_{|n}\neq \mathbf{1}
\}\right)=f_n\left(t\right) $$ and
$$\int_{0}^{t}{\mu_u\left(\{\pi \in \pinf, \pi_{|n}\neq \mathbf{1}
\}\right)du}=-\ln\left(\lambda_n\left(]t,1]\right)\right)$$
which is finite by  Hypothesis \ref{hyp}.$\Box$

\vspace*{0.5cm}

Set $\varepsilon_i=\{\{i\},\{\mcn\setminus\{i\}\}\}$ and
$\varepsilon=\sum_{i}{\delta_{\varepsilon_i}}.$ So $\varepsilon$ is a
measure on $\pinf$.  According to Bertoin
\cite{CoursBertoin03}, we know that for each  exchangeable measure
$\mu$ such that $\mu\{\un\}=0$ and $\mu\left(\{\pi \in \pinf,
\pi_{|n}\neq \mathbf{1} \}\right)<\infty$, we can find a  measure
$\nu$ on $\srb$ (dislocation measure) verifying $\nu(\un)=0$ and
$\int_{\srb}{(1-s_1)\nu(ds)}<\infty $,  and a constant $c\ge 0$
(erosion coefficient) such that
$$ \mu=\rho_{\nu}+c \varepsilon$$
where $\rho_{\nu}$ denotes the measure on $\pinf$ associated to
$\nu$ by the paint-box process.

So for $t\in [0,1[$ fixed, we can write $\mu_t=\rho_{\nu_t}+c_t
\varepsilon$ where $\nu_t$ and $c_t$ are the  instantaneous
dislocation and erosion rates of the fragmentation.

\begin{prop}
We have $\mu_t=\rho_{\nu_t}+c_t \varepsilon$ where $\nu_t$ and $c_t$ fulfill the following properties :
\begin{itemize}
\item
\begin{equation}\label{un}
\forall t\in[0,1[ \; \nu_t\left(\un\right)=0 \mbox{ and
}\int_{\srb}{\left(1-s_1\right)\nu_t\left(ds\right)}<\infty,
\end{equation}
 \item
 \begin{equation} \label{deux}
\forall u\in[0,1[ \;
\int_{0}^{u}{\int_{\srb}{\left(1-s_1\right)\nu_t\left(ds\right)dt}<\infty}
\mbox{ and }\int_{0}^{u}{c_t dt}<\infty.
\end{equation}
\end{itemize}
\end{prop}

Proof : The property (\ref{un}) is clear.  For the formula
(\ref{deux}), we shall look at the proof of the theorem in the
time-homogeneous case (cf. \cite{CoursBertoin03}). During the
proof, we obtain the following upper bound :
$$\int_{\srb}{\left(1-s_1\right)\nu_t\left(ds\right)}\le \mu_t\left(\{\pi \in \pinf, \pi_{|2}\neq \mathbf{1}
\}\right).$$ Then use Proposition \ref{intmut}.\\
 For the upper bound concerning  $c_t$ we remark :
$$c_t=\mu_t\left(\{1\},\mcn\setminus\{1\}\right)-
\rho_{\nu_t}\left(\{1\},\mcn\setminus\{1\}\right)=\mu_t\left(\{1\},\mcn\setminus\{1\}\right).\Box$$

Hence the law of a time-inhomogeneous fragmentation is
characterized by a family $(\nu_t,c_t)_{0\le t <1}$ where
$(\nu_t)_{0\le t <1}$
 and $(c_t)_{0\le t <1}$ fulfill (\ref{un}) and (\ref{deux}). One calls $\nu_t$ the instantaneous dislocation rate and $c_t$
the instantaneous erosion rate at time $t$  of the fragmentation. We will next give a probabilistic interpretation of this family.

As for the time-homogeneous fragmentations, we can construct a fragmentation with
measure $\left(\mu_t,t\in [0,1[\right)$   considering a Poisson  measure
$M$ on $  [0,1[\times\pinf\times\mcn$ with intensity
$\mu_t(d\pi)dt\otimes \sharp$ where $\sharp$ is the counting  measure.
Let $M^n$ be the restriction of $M$ to
$[0,1[\times\pn^*\times\partn$. According to Proposition
\ref{intmut}, the intensity of the  measure is finite on the
interval $[0,t]$. Then, we are in a similar case as a
time-homogeneous fragmentation
(refer to \cite{CoursBertoin03} for a proof in the homogeneous case).
Let us  rearrange the atoms of  $M^n$ according to their first coordinate.
For $n\in \mcn$, $(\pi,k)\in \pn\times \mcn$, let
$\Delta^{(.)}_n(\pi,k)$ be the following sequence of partition of $[n]$  :
$$\Delta^{(i)}_n(\pi,k)=\un \mbox{ if } i\neq k   \hspace*{0.8cm}\mbox{ and }
   \hspace*{0.8cm} \Delta^{(k)}_n(\pi,k)=\pi_{|n}.$$
We construct the   process $(\Pi_{|n}(t),t\ge 0)$
in $\pn$ with the following rules :\\
$\Pi_{|n}(0)=\un$.\\
$(\Pi_{|n}(t),t\ge 0)$ is a jump process which jumps at times
 $s$, atoms of $M^n$. More precisely, if  $(s,\pi,k)$ is an atom of $M^n$, we have
$\Pi_{|n}(s)=FRAG(\Pi_{|n}(s^-),\Delta^{(.)}_n(\pi,k))$.
We can then check that this construction is compatible with the restriction and the constructed process is a fragmentation with
 measure
$\mu_t$.

We have also a Poissonian construction for a mass-fragmentation (cf. \cite{Berestycki02}).
First we use that if $F=(F(t),t\in [0,1[)$ is a mass-fragmentation with parameters $(\nu_t,0)_{0\le t<1}$, then $\tilde{F}=(e^{-\int_{0}^{t}{c_s ds}}F(t),t\in [0,1[)$ is
a mass-fragmentation with parameters $(\nu_t,c_t)_{0\le t <1}$. So, we  remark that the family  of instantaneous erosion coefficients plays only a deterministic role in
the fragmentation.
To find a Poissonian construction for the  mass-fragmentations $(F(t),t\in [0,1[)$ with parameters $(\nu_t,0)_{0\le t <1}$,
consider then a fragmentation on partitions $(\Pi(t),t\in [0,1[)$ such that $F=\Lambda(\Pi)$ where
 $\Lambda$ is the application which associates to a partition its frequency sequence. So $\Pi$ can be constructed from a Poisson measure $M$.
 Consider $K$, image of $M$ by the
application
\begin{eqnarray*}
\pinf \times \mcn & \longrightarrow &\sr \times \mcn\cup \infty\\
(\Delta(\cdot),k(\cdot))&\longmapsto &
(\Lambda(\Delta(\cdot)),f(\cdot,k(\cdot))),
\end{eqnarray*}
where   $f$ is the  function which associates to $k$ the frequency rank of the
block $B_k(t^{-})$. \\
Berestycki  \cite{Berestycki02} then  proves that   $K$ is a Poisson measure on $[0,1[\times
\sr \times \mcn$ with intensity measure  $\nu(ds)dt\otimes \sharp$. \\
Set
$$K=(t,S(t),k(t))_{t\in [0,1[}=(t,(s_1(t),s_2(t),\ldots),k(t))_{t\in [0,1[}.$$
Then, if $(t,S(t),k(t))$ is an atom of $K$, then at time $t$,
the $k(t)$-th largest block of the  fragmentation at time
$t^{-}$ will be fragmented according to $S(t)$.

\vspace*{0.5cm}

Let us  now determine the effects of a deterministic change-time on a fragmentation.

\begin{prop}\label{changtemps}
Let  $(\Pi(t),t\in [0,1[)$ be a fragmentation with parameter
$(c_t,\nu_t)$. Set $\Pi'(t)=\Pi(\beta(t))$ where
$\beta:[0,1[\rightarrow \Rp$ is a  strictly
increasing derivable function. Let $J$ be the image of $[0,1[$ by $\beta$
($J$ is thus an interval of $\Rp$). \\
Then $(\Pi'(t),t\in J)$ is a fragmentation with
parameter $(c'_t;\nu'_t)_{t\in J}$ where
$$c'_t=\beta'(t)c_{\beta(t)} \hspace*{2cm}
\nu'_t=\beta'(t)\nu_{\beta(t)}.$$
\end{prop}

Proof. A Markov process remains a Markov process after
a deterministic time-change. The law of
$\Pi'(t+t')$ given $\Pi'(t)=\pi$, is
$FRAG(\pi,\pi^{(\cdot)})$, where $\pi^{(\cdot)}$ is an iid sequence with
law $\mathbb{P}_{\beta(t),\beta(t+t')-\beta(t)}$. Thus $\Pi'$ is a fragmentation.\\
Let us  calculate its jump rates $q'_{\pi,t}$.
\begin{eqnarray*}
q'_{\pi,t}dt&=&\mcp\left(\Pi'_{|n}\left(\tau'_n\right)=\pi \; \&
\;\tau'_n \in [t,t+dt]\;|\; \tau'_n\ge t\right)\\
&=&\mcp\left(\Pi_{|n}\left(\beta(\tau'_n)\right)=\pi \; \&
\;\tau'_n \in [t,t+dt]\;|\; \tau'_n\ge t\right)\\
&=&\mcp\left(\Pi_{|n}\left(\tau_n\right)=\pi \; \& \;\tau_n \in
[\beta(t),\beta(t+dt)]\;|\; \tau_n\ge \beta(t)\right)\\
&\sim&\mcp\left(\Pi_{|n}\left(\tau_n\right)=\pi \; \& \;\tau_n \in
[\beta(t),\beta(t)+\beta'(t)dt)]\;|\; \tau_n\ge \beta(t)\right)\\
&\sim& \beta'(t)q_{\pi,\beta(t)}dt.
\end{eqnarray*}

So $q'_{\pi,t}= \beta'(t)q_{\pi,\beta(t)}$. We thus deduce similar
relations between
 $\nu_t$ and $\nu'_t$ and between $c_t$
and $c'_t$. $\Box$

\subsection{Law of the tagged fragment}\label{tagged}

An application of the above decomposition is for example to calculate the law of the frequency  of the  block containing $1$,
$|\Pi_1\left(t\right)|$, for an exchangeable standard fragmentation.
 This quantity is interesting because it represents the law of a  size-biased picked block.
 We have the following theorem :

\begin{theo}\label{sub}
There exists a process $(\xi(t),t\in [0,1[)$ with independent
increments such that $|\Pi_1\left(t\right)|=
\exp\left(-\xi_t\right)$. Its law is characterized by the identity :
$$
\mce\Big(|\Pi_1\left(t\right)|^q\Big)=\mce\Big(\exp (-q\xi_t)\Big)=\exp\left(-\int_{0}^{t}{\phi_u\left(q\right)du}\right), \hspace*{1cm} q>0$$
$$\mbox{where }\phi_t\left(q\right)=c_t\left(q+1\right)+\int_{\sr}{\left(1-\sum_{i=1}^{\infty}
{s_i^{q+1}}\right)\nu_t\left(ds\right)}.$$
\end{theo}

In the sequel, we will  also use the notation
$\psi(t,q)=\int_{0}^{t}{\phi_u\left(q\right)du}.$

\vspace*{0.3cm}

This result is very close to the corresponding result in the
homogeneous case. We just loose the stationarity of the increments
of $\xi(t)$. The demonstration itself is similar to the
homogeneous case and we just sketch the proof here. For more
details, refer to \cite{CoursBertoin03}.

We use the equality :
$$\mcp[\Pi_{|k+1}\left(t\right)=\mathbf{1}]=\mce[|\Pi_1\left(t\right)|^k],$$
which we get by conditioning on $|\Pi_1\left(t\right)|$. Then
remark the event $\{\Pi_{|k+1}\left(t\right)=\mathbf{1}\}$
corresponds, looking at the Poissonian construction, to an absence
of Poisson atom in the subset $[0,t]\times\{\pi \in
\pinf,\pi_{|k+1}\left(t\right)\neq\mathbf{1}\}\times\{1\}$.\\
So the formula is true for every positive integer. Besides, we
remark that the law of $|\Pi_1\left(t\right)|$ is characterized by
its moments, thanks to the independence of the increments (when you take
the logarithm)
 and because the process takes values in $[0,1]$.\\
 By uniqueness of the analytic continuation, we deduce that the formula is true for every $q>0$.
And by the monotone convergence theorem,
$\psi(t,q)$ is continuous in $q$ at 0. $\Box$

Thanks to this formula, we can characterize the processes which
have proper frequencies, i.e. with
$\sum_{i=1}^{\infty}{|\pi_i|}=1$.

\begin{prop}\label{propre}
We have :
$$\;\mcp\left(\Pi\left(t\right) \mbox{ is proper }\right)=1\Leftrightarrow
\left(c_u=0 \mbox{ and } \nu_u\left(\sum_{i}{s_i}<1\right)=0
\mbox{ for }0\le u\le t \;a.e.\right).$$
\end{prop}

Proof. First remark
$$\lim_{k\rightarrow 0}\mce[|\Pi_1\left(t\right)|^k]=\lim_{k\rightarrow 0}\mce[|\Pi_1\left(t\right)|^k \mathbf{1}_{|\Pi_1\left(t\right)|\neq
0}]=\mce[\mathbf{1}_{|\Pi_1\left(t\right)|\neq
0}]=1-\mcp\left(|\Pi_1\left(t\right)|=0\right).$$

Then  we have :
\begin{eqnarray*}
\mcp\left(\Pi\left(t\right) \mbox{ is  proper }\right)=1 &
\Leftrightarrow &
\mcp\left(|\Pi_1\left(t\right)|=0\right)=0\\
& \Leftrightarrow & \exp\left(-\psi(t,0)\right)=1\\
& \Leftrightarrow &\psi(t,0)=0\\
 & \Leftrightarrow &  \phi_u\left(0\right)=0 \mbox{
for }0\le u\le t \;a.e. \Box
\end{eqnarray*}

Recall from \cite{Bertoin03} that if  $(X\left(t\right),t\in [0,1[)$ is a
time-homogeneous mass-fragmentation,  $\phi$  the  Laplace
exponent associated to the tagged fragment and  $\mathcal{F}_t=\sigma\left(X\left(s\right),s\le
t\right)$, then
$$\exp\left(t\phi\left(p\right)\right)\sum_{i=1}^{\infty}{X_i^{p+1}\left(t\right)}\mbox{ is
a }\mathcal{F}_t\mbox{-martingale}.$$

We can obtain a similar theorem in the time-inhomogeneous case.

\begin{prop}Consider  $(\Pi(t),t\in [0,1[)$ a time-inhomogeneous fragmentation on partitions. Let
$X\left(t\right)=\left(X_i\left(t\right)\right) \in \srb$ be its decreasing sequence of frequencies.
Set $\mathcal{F}_t=\sigma\left(X\left(u\right);u\le t\right)$. Let
$\phi_u$ be its instantaneous Laplace exponent and $\psi(t,p)=\int_{0}^{t}\phi_u(p)du$. Then
$$M(t,p)=\exp\left(\psi(t,p)\right)\sum_{i=1}^{\infty}{X_i^{p+1}\left(t\right)}\mbox{ is
a } \mathcal{F}_t\mbox{-martingale.}$$
\end{prop}

Proof. It is the same idea as in the time-homogeneous case. Set $\mathcal{G}_t=\sigma\left(\Pi\left(u\right),u\le t\right)$.
Then
  $\mathcal{E}\left(t,p\right)=\exp(-p\xi_t+\psi(t,p))$ is an  $\mathcal{G}_t$-martingale
and we  remark that $M(t,p)$ is the projection of $\mathcal{E}\left(t,p\right)$ on $\mathcal{F}_t$. $\Box$

\section{Application to Ruelle's cascades }\label{appl}

\subsection{Jump rates of Ruelle's fragmentation}

Let $(\Pi(t),t\in [0,1[)$ be Ruelle's fragmentation with values in
partitions. For each integer $n$, $(\Pi_{|n},t\in[0,1[)$ is a Markov process in the finite
space of partition of $[n]$. The law of such a process is entirely determined by its
jump rates from one state to another.

Let us  calculate its jump rates.
Set
$\pi=\left(\pi_1,\ldots,\pi_k\right)\in\pn^*$.  Fix $t\in [0,1]$. Let
 $q_{ t}\left(n_1,\ldots,n_k\right)$ be the probability that
$\Pi_{|n}\left(t\right)$ has blocks with size
$\left(n_1,\ldots,n_k\right)$. Recall (cf. Proposition \ref{EPFBS}) that
$$ q_{
t}\left(n_1,\ldots,n_k\right)=\frac{\left(k-1\right)!}{\left(n-1\right)!}t^{k-1}\prod_{i=1}^{k}{[1-t]_{n_i-1}}.$$
So from Proposition \ref{frag} and \ref{alpha-theta}
\begin{eqnarray*}
\mcp\left(\tau_n \in [t,t+s],
\pic_{|n}\left(\tau_n\right)=\pi\;|\; \tau_n\ge t\right)&=&
\mcp\left(\pic_{|n}\left(t+s\right)=\pi  \;|\;
\pic_{|n}\left(t\right)=\mathbf{1}\right)\\
&=&p_{t+s,-t}(n_1,\ldots,n_k)\\
&=& \frac{\left[\frac{-t}{t+s}\right]_k}{[-t]_n}\prod_{i=1}^{k}{-[-t]_{n_i}}\\
&\sim& s\frac{(-1)^{k+1}(k-2)!\prod_{i=1}^{k}{[-t]_{n_i}}}{t
[-t]_{n}}.
 \end{eqnarray*}

Remark that we could also have calculated this quantity using Proposition \ref{fragcog} and Bayes' Formula.
So we obtain the following proposition :
 \begin{prop}\label{qpit}
For $\pi=\left(\pi_1,\ldots,\pi_k\right)\in\pn^*$ and for $t \in
[0,1[$ we have :
$$q_{\pi,t}=\frac{q_{t}\left(n_1,\ldots,n_k\right)}{t\left(k-1\right)q_{t}\left(n\right)}.$$
 \end{prop}

 \subsection{Instantaneous erosion coefficient and dislocation measure}

It is well known that the Bolthausen-Sznitman's coalescent is a
process with proper frequencies (cf. Proposition \ref{propre}).
So, the erosion coefficient $c_t$ should be identically
 zero. We can check this with a short calculation.
In fact, consider
$\pi=\varepsilon_1=\Big\{\{1\},\mcn\setminus\{1\}\Big\}$ and
$\pi_n=\pi_{|n}.$\\
According to Proposition \ref{qpit}, we have $q_{\pi_n,t}=\frac{
q_{t}\left(1,n-1\right)}{t q_{t}\left(n\right)}=\frac{1}{ n-1-t}$.
And
$c_t=\lim_{n\rightarrow \infty} q_{\pi_n,t}=0$. \\
Thus $c_t=0$ for all $t \in [0,1[$.

\vspace*{0.5cm}

Let us denote by $\tilde{\sr}$ the set of  the positive sequence
with sum $1$. From a measure
$\eta$ on $\tilde{\sr}$, we can define a measure $p$ on $\pinf$ (cf. \cite{Pitman02} p. 61) :\\
Conditionally on a sequence $(s_i,i\ge 1)$ drawn with respect to the measure $\eta$, we construct the following law on partitions :\\
$1$ is in the first block. Fix $n\ge 1$.  Suppose $\Pi_n$ has $k$
blocks. The integer $n+1$ will be :
\begin{itemize}
\item in the block $j$ with probability $s_j$ (for $j\le k$),
\item in a new  block with probability $1-\sum_{i=1}^{k}{s_i}$.
\end{itemize}

So we have

\begin{equation}\label{g}
p(\pi)=\mce^{\eta}\Big(\prod_{i=1}^{k}{s_i^{n_i-1}}\prod_{i=1}^{k-1}
{(1-\sum_{j=1}^{i}{s_j})}\Big),
\end{equation}
where $\pi=(\pi_1,\ldots,\pi_k)$ et $|\pi_i|=n_i$.

If the measure
 $\eta$ is a dislocation measure (i.e verifies $\int_{\sr}{(1-s_1)\eta(ds)}<\infty$),
 then $p$ is finite on $\pn^{*}$.
In fact, for all $k\ge 2$, we have \\
$\prod_{i=1}^{k}{s_i^{n_i-1}}\prod_{i=1}^{k-1}
{(1-\sum_{j=1}^{i}{s_j})}\le 1-s_1$.

Let us now look at the dislocation measure. In this direction, let us
introduce the following measure :

\begin{defi}
Fix $\alpha \in ]0,1[$. Consider the measure $\eta_{\alpha}$
defined as follows on $\tilde{S}$ : first,
$$\eta_{\alpha}(s_1\in dx)=\alpha
 x^{-\alpha}(1-x)^{-1}\un_{0<x<1}dx,$$
 and second, conditionally  on $s_1=x$, the sequence $(s_{i+1}/(1-x),i\in
 \mcn)$ has the law of a random variable with law
  $PD(\alpha,0)$ of which the terms have been size-biased rearranged.
 We denote $PD(\alpha,-\alpha)$ the image of $\eta_\alpha$ by ranking the $s_i$ in the decreasing
 order. $PD(\alpha,-\alpha)$ is then an infinite measure on
 $\sr$.
 \end{defi}

Remark that
the construction of the measure
 $PD(\alpha,-\alpha)$ is similar, except for the
  normalization, to the construction of a Poisson-Dirichlet
measure with the forbidden parameter $\theta=-\alpha$.

\begin{prop}
 Define $p_\alpha$ as the measure on $\pinf$ associated to
 $\eta_\alpha$ as above. Then
$p_\alpha$ is an exchangeable measure on $\pinf$. Its EPPF for the
partitions non-reduced to one block is :
\begin{equation}\label{EPMF}
p_\alpha(n_1,\ldots,n_k)=\frac{(k-2)!}{-[-\alpha]_{n}}\prod_{i=1}^{k}{-[-\alpha]_{n_i}} \mbox { for all } k\ge 2.
\end{equation}
\end{prop}

Proof. Let us  first check $\int_{\sr}{(1-s_1)\eta_\alpha(ds)<\infty}$.
\begin{equation}\label{maj}
\int_{\sr}{(1-s_1)\eta_\alpha(ds)}=\int_{0}^{1}{(1-s_1)\alpha s_1^{-\alpha}(1-s_1)^{-1}ds_1}
=\frac{\alpha}{1-\alpha}.
\end{equation}
Using formula (\ref{g}) and  the definition of
$\eta_\alpha$, we have :
\begin{eqnarray*}
p_\alpha(\pi)&=&\left(\int_{0}^{1}{x^{n_1-1}(1-x)^{\sum_{i=2}^{k}{n_i}}\eta_\alpha(s_1
\in dx)}\right)p_{\alpha,0}(n_2,\ldots,n_k)\\
&=&\alpha
\left(\int_{0}^{1}{x^{n_1-1-\alpha}(1-x)^{n-n_1-1}dx}\right)p_{\alpha,0}(n_2,\ldots,n_k)\\
&=& \alpha
\frac{\Gamma(n_1-\alpha)\Gamma(n-n_1)}{\Gamma(n-\alpha)}\frac{(k-2)!}{\alpha
(n-n_1-1)!}\prod_{i=2}^{k}{-[-\alpha]_{n_i}}\hspace*{1cm} \mbox{ according to (\ref{alpheq})} \\
&=& \frac{[-\alpha]_{n_1}}{[-\alpha]_{n}}(k-2)!
\prod_{i=2}^{k}{-[-\alpha]_{n_i}} \\
&=&\frac{(k-2)!}{-[-\alpha]_{n}}\prod_{i=1}^{k}{-[-\alpha]_{n_i}}.
\end{eqnarray*}

So, we find the foretold formula and this one is symmetric in the
variables
 $(n_1,\ldots,n_k)$, thus the measure is an exchangeable measure
(cf. \cite{Pitman02} Theorem 24). We also deduce that
$\eta_{\alpha}$ is the image of $PD(\alpha,-\alpha)$ by a
size-biaised reordering and $p_\alpha=\rho_{PD(\alpha,-\alpha)}$
(where $\rho_{PD(\alpha,-\alpha)}$ is the measure on $\pinf$
obtained from $PD(\alpha,-\alpha)$ by the paint-box
construction.)$\Box$

Next, we observe that for every partition  $\pi$ not reduced to
one block, we have
$$q_{\pi,t}=\frac{1}{t}p_{t}(\pi).$$

Indeed, this follows from Proposition
 \ref{qpit} and  formula (\ref{EPFBSF}) of Pitman.
In conclusion, we may now state the following theorem :

\begin{theo}\label{mesure}
The instantaneous dislocation measure
 $\nu_t$ of Ruelle's fragmentation at
time $t$ is given by :
$$\nu_t=\frac{1}{t}PD(t,-t).$$
\end{theo}

\subsection{Absolute continuity of the dislocation  measure and  $PD(\alpha,0)$}

Let us recall that, if $\Pi$ is a random partition with law
$p_{\alpha,0}$ and   $K_n$  the number of block of $\Pi_{|n}$, then the
limit of $K_n/n^\alpha$ exists almost surely and has the
Mittag-Leffler law with index $\alpha$ (cf. \cite{Pitman02}
Theorem 31)

\begin{prop}
For each $\alpha \in ]0,1[$ the  measure $p_\alpha$ is absolutely
continuous with respect to the measure $p_{\alpha,0}$. More
precisely, we have :
$$p_\alpha(d\pi)=\Gamma(1-\alpha)S_{\alpha}^{-1}p_{\alpha,0}(d\pi)
\hspace*{0.8cm}\mbox{ where }S_{\alpha}=\underset{n \rightarrow
\infty}{\lim}\frac{K_n}{n^\alpha}.$$
\end{prop}

Proof. Let $(\mathcal{F}_{n})_{n\ge1}$ be the  filtration of $\Pi_{|n}$.\\
Fix $k\ge2$. Set $p^k_{\alpha}= p_{\alpha}\un_{\mathcal{P}_k^*}$.
We consider
$$M^k_{\alpha,n}={\frac{dp^k_{\alpha}}{dp_{\alpha,0}}}{\Big|\mathcal{F}_n}.$$
Using formula (\ref{EPFBSF}) and  (\ref{EPMF}), we have :
$$M^k_{\alpha,n}=\frac{\Gamma(1-\alpha)\Gamma(n)}{\Gamma(n-\alpha)(K_n-1)}
\un_{\mathcal{P}_k^*}\hspace*{1cm}\mbox{ for } n\ge k,$$ where
$K_n$ denotes the number of block of $\Pi_{|n}$.
$M^k_{\alpha,n}$ is a positive martingale, thus it converges almost surely to a random variable  $M^k_{\alpha}$.\\
Let now use $$\frac{K_n}{n^\alpha}\rightarrow
S_{\alpha}\hspace*{0.2cm} \mcp_{\alpha,0}-\mbox{a.s.}
\hspace*{0.3cm} \mbox{ and }\hspace*{0.3cm}
\frac{\Gamma(1-\alpha)\Gamma(n)}{\Gamma(n-\alpha)(K_n-1)}\sim
\frac{\Gamma(1-\alpha)n^\alpha}{K_n}.$$ We deduce
$$M^k_\alpha=\frac{dp^k_{\alpha}}{dp_{\alpha,0}}=\Gamma(1-\alpha)S_{\alpha}^{-1}\un_{\mathcal{P}_k^*}\hspace*{0.2cm}
\mcp_{\alpha,0}-\mbox{a.s.}$$

So, according to  martingale theory (cf. \cite{Durrett91}
p.210),  for all $A\subset \mathcal{P}_k^*$, we have  :

$$p_{\alpha}(A)=\mce_{\alpha,0}\left(\Gamma(1-\alpha)S_\alpha^{-1}\un_{A} \right)
+p_{\alpha}(A\cap \{S=0\}),$$
where $S=\limsup \frac{K_n}{n^\alpha}$.\\
Set $x \in ]0,1[$. Let us define $q_{\alpha}(\cdot)=c
p_{\alpha}(\;\cdot\;|\;|\Pi_1|=x)$ where  $c$ is chosen such that
$q_\alpha$ is a probability. Let $s=(s_1,\ldots)\in \srb$ be the
frequency sequence of a partition with law $q_\alpha$. According
to the construction of $p_\alpha$, we have
$$(s_{i+1})_{i\in \mcn}\overset{law}{=}(1-x)(p_i)_{i\in \mcn}, $$
where $(p_i)_{i\in \mcn}$ has the  $PD(\alpha,0)$ law.

According to Lemma 34 of Pitman's course \cite{Pitman02}, for a
random partition, $S$ exists and belongs almost surely to
$]0,\infty[$ iff there exists $Z$ random variable on $]0,\infty[$
such that  $P_i\sim Z i^{-1/\alpha}$,
 where $P_i$ is the decreasing sequence of the frequencies.
Here we know the existence of such a random variable $Z \in ]0,\infty[$ for a $PD(\alpha,0)$ law. Set $Y=(1-x)Z$ then
$$s_i\sim Y i^{-1/\alpha}.$$
So we have
$$p_{\alpha}(S=0\;|\;|\Pi_1|=x)=0.$$
Thus $$p_{\alpha}(S=0)=0.$$ We conclude that $$\forall A \in \pinf \mbox{ such
that } \un \not\in \overline{A}\hspace*{0.5cm} p_{\alpha}(A)=
\mce_{\alpha,0}\left(\Gamma(1-\alpha)S_\alpha^{-1}\un_A \right).
\Box$$

\begin{theo}
The dislocation measure of Ruelle's fragmentation at time $t$ is
 absolutely continuous with respect to
the measure $PD(t,0)$. More precisely, we have for all continuous function  $f$ on $\srb$ :
$$\nu_t(f)=\frac{1}{t}\mce_{(t,0)}\left(L_{t}^{-1}f(V)\right)$$
where $L_{\alpha}=\underset{n\rightarrow \infty}{\lim}{nV_n^\alpha}$.
\end{theo}

Proof. We use that if $(s_i)_{i\ge 1} \in \sr$ is the frequency sequence of an $(\alpha,0)$-partition $\Pi_{\infty}$,
then $\Gamma(1-\alpha)L_{\alpha}$ exists almost surely and it is equal almost surely to
 $S_\alpha=\underset{n\rightarrow \infty }{\lim}{\frac{K_n}{n^\alpha}}$ (cf. \cite{Pitman02} Theorem 36).\\
Use Theorem \ref{mesure} to finish the proof. $\Box$

\begin{rem}
$L_\alpha$ is not a continuous function on $\sr$.
\end{rem}

\subsection{Law of the tagged  fragment}
In this section, we determine the law of the tagged fragment. Actually, its law has already been determined by Pitman \cite{Pitman99}. He proves that
 $|\Pi_1\left(t\right)|$ has a
$\beta\left(1-t,t\right)$  law.
 So we  check that we find the same result.

Hence, according to Section \ref{tagged}, we shall calculate
$\phi_t\left(k\right)=\int_{\sr}{\left(1-\sum_{i=1}^{\infty}{s_i^{k+1}}\right)\nu_t\left(ds\right)}$.
Recall  that $p_t$ denotes the measure on $\pinf$
associated to the measure $PD(t,-t)$.\\
We have
$$
\mce[|\Pi_1\left(t\right)|^k]=\exp\left(-\int_{0}^{t}{\phi_u\left(k\right)du}\right).$$

Thus
\begin{eqnarray*}
\phi_t\left(k\right)&=&\mce_{v_t}\left(\rho_s\left(\Pi_{|k+1}\neq \un\right)\;|s\right)\\
&=&\frac{1}{t}p_t\left(\Pi_{|k+1}\neq \un\right).
\end{eqnarray*}

So we must calculate  $p_t\left(\Pi_{|k+1}\neq \un\right)$. We
will do this recursively.\\
For  $k=1$, we have
$$p_t\left(\Pi_{|2}\neq
\un\right)=\frac{[-t]_1^2}{-[-t]_2}=\frac{t}{1-t},$$ and for $k\ge
2$ $$p_t\left(\Pi_{|k+1}\neq \un\right)=p_t\left(\Pi_{|k}\neq
\un\right)+p_t\left(\Pi_{|k+1}=\Big\{\{1,\ldots,k\},\{k+1\}\Big\}\right)=p_t\left(\Pi_{|k}\neq
\un\right)+\frac{t}{k-t}.$$

Thus we have :
$$ p_t\left(\Pi_{|k+1}\neq \un\right)=\sum_{i=1}^k
\frac{t}{i-t} \hspace*{0.5cm}\mbox{ and so }\hspace*{0.5cm}
\int_{0}^{t}{\phi_u\left(k\right)du}=
\ln\left(\prod_{i=1}^{k}{\frac{i}{i-t}}\right).$$

So we deduce

$$\mce[|\Pi_1\left(t\right)|^k]=\prod_{i=1}^{k}{\frac{i-t}{i}}.$$
The right-hand side coincides with the $k$-th moment of a
$\beta\left(1-t,t\right)$ law. So $|\Pi_1\left(t\right)|$ has a
$\beta\left(1-t,t\right)$ law and we deduce :

$$\forall k>0,\;\mce[|\Pi_1\left(t\right)|^k]=\frac{\Gamma\left(k+1-t\right)}{\Gamma\left(1-t\right)
\Gamma\left(k+1\right)}.$$

More generally, we can determine the law of the process
$\left(|\Pi_1\left(t\right)|,t \in [0,1[\right)$. By the
homogeneous property of fragmentation in space, the process
$\left(\frac{|\Pi_1\left(t+s\right)|}{|\Pi_1\left(t\right)|},s \in
[0,1-t[\right)$ is independent of $|\Pi_1\left(t\right)|$ (cf.
Theorem \ref{sub}). So we can calculate the finite dimensional law
of the process $\left(|\Pi_1\left(t\right)|,t \in [0,1[ \right)$
and we deduce that the process has the same law as the process
$\left(\frac{\gamma\left(1-t\right)}{\gamma\left(1\right)},t \in
[0,1[\right)$ (result already proved by Pitman \cite{Pitman99}).

\begin{rem}
We have also an expression for $\psi(t,p)$ :
$$\psi(t,p)=\ln \left( \frac{\Gamma\left(1-t\right)
\Gamma\left(k+1\right)}{\Gamma\left(k+1-t\right)}\right).$$
\end{rem}

\section{Behavior of the fragmentation at large and small times}\label{asymp}

\subsection{Convergence of the empirical measure }
Let
 $(\Pi(t),t\in[0,1[)$ be a Ruelle's fragmentation on the partitions.
  Let $(X(t),t\in[0,1[)$, $X(t)=(X_i(t))_{i\ge 1}\in\sr$ be its process of ranked frequencies.\\
We are interested in the empirical measure  $\rho_t$ defined by
:

$$ \rho_t=\sum_{i=1}^{\infty}{X_i(t)\delta_{(t-1) \ln X_i(t)}}.$$

\begin{prop}
For every bounded continuous function  $f$ on $\Rp$:
$$\underset{t \rightarrow
1}{\lim}{\int f(y)\rho_t(dy)}=\int_{0}^{\infty}{f(y)e^{-y}dy}\;
\mbox{ in } L^{2}.$$
\end{prop}

We split the proof in two parts. We will successively prove the
 following two points:
\begin{equation}\label{l1}
 \underset{t \rightarrow 1}{\lim}{\mce\left(\int
f(y)\rho_t(dy)\right)}=\int_{0}^{\infty}{f(y)e^{-y}dy},
\end{equation}
\begin{equation}\label{l2}
\underset{t \rightarrow 1}{\lim}{\mce\left[\left(\int
f(y)\rho_t(dy)\right)^2\right]}=\left(\int_{0}^{\infty}{f(y)e^{-y}dy}\right)^2.
\end{equation}

Set $\xi_t=-\ln |\Pi_1(t)|$.
Let us  recall
$$|\Pi_1\left(t\right)|\sim \beta\left(1-t,t\right),$$
and observe  :
$$\mce\left(\int
f(y)\rho_t(dy)\right)=\mce\Big( f((1-t) \xi_t)\Big).$$

The following lemma  clearly implies (\ref{l1}).
\begin{duge}
Set $\xi_t=-\ln |\Pi_1(t)|$ where $\Pi(t)$ is the Ruelle's fragmentation. Then
$$\underset{t \rightarrow
1}{\lim}{(1- t )\xi_t}=\mathbf{e} \mbox{ in distribution }$$
 where $\mathbf{e}$ denotes the exponential law with parameter $1$.
\end{duge}

Proof.
Let us calculate the Laplace transform  of $(1-t) \xi_t$.
\begin{eqnarray*}
\mce\left(e^{-q(1- t) \xi_t}\right)&=&\mce\left(|\Pi_1(t)|^{q(1-t)}\right)\\
&=&\frac{\Gamma\left(q(1-t) +1-t\right)}{\Gamma\left(1-t\right)
\Gamma\left(q(1-t)+1\right)}\\
&\underset{t \rightarrow 1}{\longrightarrow}& \frac{1}{q+1}.
\end{eqnarray*}

Since $\frac{1}{q+1}$ is the Laplace transform of the exponential
law, by  Lévy's Theorem, $(1-t) \xi_t$ converges in law to
$\mathbf{e}$. $\Box$

To prove (\ref{l2}), we consider $\xi^{'}_t=-\ln|\Pi_2(t)|$ where
$\Pi_2(t)$ is the block containing the integer $2$. Observe that $\xi_t$ and
$\xi^{'}_t$ have the same law but are not independent, and that
$$\mce\left[\left(\int f(y)\rho_t(dy)\right)^2\right]=\mce\left[
f\Big((1- t)\xi_t\Big)f\left((1- t)
\xi^{'}_t\right)\right].$$

Set $T=\inf\left\{t>0, \Pi_1(t)\neq\Pi_2(t)\right\}$, so   $T$ is
 almost surely finite and  conditionally on $T$,  $\xi_T$ and
$\xi^{'}_T$, the processes $(\xi_t,t\ge T)$ and $(\xi^{'}_t,t\ge
T)$ are independent. From this, we deduce (\ref{l2}) and then the
$L^{2}$-convergence of $\int{f(y)\rho_t(dy)}$ (refer to
\cite{Bertoin03} for details).$\Box$

So, informally, this proposition proves that, if we consider the size of a typical fragment $X(t)$,  then, as $t$ tends to $1$, we have
$$|\log X(t)|\sim \frac{C}{1-t}$$
where $C$ is a random factor.

\subsection{Additive martingale }

In this section, we aim at studying the convergence of the
martingale $M(t,p)$ defined in Section \ref{tagged} and we follow the ideas of Bertoin and Rouault
\cite{Bertoinrouault03} who introduce a new probability to
prove the convergence.

Recall the  following notation :\\
 $\mathcal{F}_t=\sigma\left(X_i\left(u\right),u\le t\right)$ is the filtration
of the frequency sequence.\\
$\mathcal{G}_t=\sigma\left(\Pi\left(u\right),u\le t\right)$ is the
filtration of the fragmentation process on the partitions.\\
So we have $\mathcal{F}_t\subseteq \mathcal{G}_t$.\\
Set $\xi_t=-\ln\left(|\Pi_1\left(t\right)|\right)$. It is an increasing process
with independent increments.\\
$M\left(t,p\right)=\exp\left(\psi\left(t,p\right)\right)\sum_{i=1}^{\infty}{|X_i\left(t\right)|^{p+1}}.$
$M\left(\cdot,p\right)$ is then a $\mathcal{F}_t$-martingale. \\
$\mathcal{E}\left(t,p\right)=\exp\left(\psi\left(t,p\right)-p
\xi_t\right).$
$\mathcal{E}\left(\cdot,p\right)$ is  a $\mathcal{G}_t$-martingale. \\
As $\mce(|\Pi_1(t)|^p\;|\;X(t))=\sum_{i}{X_i(t)^{p+1}}$, we have
 $\mce\left(\mathcal{E}\left(t,p\right) \;|\;
\mathcal{F}_t\right)=M\left(t,p\right)$.\\
We denote  $\mcq$ the  probability on  $\mathcal{G}$ defined by :
$$d\mcq_{|\mathcal{G}_t}=\mathcal{E}\left(t,p\right)d\mcp_{|\mathcal{G}_t}. \mbox{ So we have also }
d\mcq_{|\mathcal{F}_t}=M\left(t,p\right)d\mcp_{|\mathcal{F}_t}. $$

\begin{prop}
Fix $p>0$. We have :
$$\lim_{t\rightarrow 1}{M(t,p)}=0 \hspace*{0.3cm}\mcp\mbox{-a.s.}$$
\end{prop}

Proof. A martingale theorem (cf. \cite{Durrett91} p.210) asserts
that  if $\limsup M\left(t,p\right)=\infty \;\mcq$-a.s., then
$\lim M\left(t,p\right)=0 \;\mcp$-a.s.

We have  $$M\left(t,p\right)\ge
\exp\left(\psi\left(t,p\right)\right)|\Pi_1\left(t\right)|^{p+1}=\exp\left(\psi\left(t,p\right)-\left(p+1\right)\xi_t\right).$$
Set $N_t=\psi\left(t,p\right)-\left(p+1\right)\xi_t$. We will
prove that $\limsup N_t=\infty \;\mcq$-a.s.

Let us recall that, under $\mcp$, $|\Pi_1\left(t\right)|$ has
$\beta\left(1-t,t\right)$ law. So for all  $\lambda\ge 0$ we have :

$$\mcq\left(\xi_t\ge \lambda\right)=
\mce^{\mcp}\left(\mathcal{E}\left(t,p\right)\un_{\{\xi_t\ge
\lambda\}}\right)=
 \frac{\Gamma\left(p+1\right)}{\Gamma\left(p+1-t\right)\Gamma\left(t\right)}
 \int_{0}^{e^{-\lambda}}{x^{p-t}\left(1-x\right)^{t-1}dx}.
$$
So for $A\le \psi\left(t,p\right),$
$$\mcq\left(N_t\le A\right)=\mcq\left(\xi_t\ge \frac{\psi\left(t,p\right)-A}{p+1}\right)=
\frac{\Gamma\left(p+1\right)}{\Gamma\left(p+1-t\right)\Gamma\left(t\right)}
\int_{0}^{e^{-\frac{\psi\left(t,p\right)-A}{p+1}}}
{x^{p-t}\left(1-x\right)^{t-1}dx}.
$$

Recall $\psi\left(t,p\right)\sim -\ln (1-t)$ as $t\uparrow 1$.
Choose $A\left(t\right)=-\frac{1}{3}\ln (1-t)$. So for
$t$ large enough, we have $ \psi\left(t,p\right)-A\left(t\right)\ge -\frac{1}{3}\ln (1-t)$. \\
Set $g\left(t\right)=(1-t)^{\frac{1}{3\left(p+1\right)}}$. We have :
\begin{eqnarray*}
\mcq\left(N_t\le A\left(t\right)\right)&\le&
\frac{\Gamma\left(p+1\right)}{\Gamma\left(p+1-t\right)\Gamma\left(t\right)}\int_{0}^{g\left(t\right)}
{x^{p-t}\left(1-x\right)^{t-1}dx}\\
&\le&
\frac{\Gamma\left(p+1\right)}{\Gamma\left(p+1-t\right)\Gamma\left(t\right)}\left(1-g\left(t\right)\right)^{t-1}\frac{1}{p+1-t}g\left(t\right)^{p+1-t}\\
&\le& \varepsilon_p\left(t\right),
\end{eqnarray*}
where $\varepsilon_p\left(t\right)$ is a function with limit 0 at
$t=1$.

So $\underset{t \rightarrow 1}{\lim} \mcq\left(N_t\ge
A\left(t\right)\right)=1$ and then $\mcq\left(\limsup N_t
<\infty\right)=0$. We deduce :
$$\limsup_{t\rightarrow 1} M\left(p,t\right)=\limsup_{t\rightarrow 1} N\left(p,t\right) =\infty \; \mcq\mbox{-a.s. } \mbox{ and so
} \lim_{t\rightarrow 1} M\left(p,t\right)=0 \; \mcp\mbox{-a.s. } \Box$$

\begin{rem}
In the case $p=0$, as the process has proper frequencies, we have
$M\left(0,t\right)=1 \; \mcp$-as for all $t\in [0,1[$.
\end{rem}

\subsection{Small times behavior}

In this section, we obtain information on the behavior of the two
largest blocks of  Ruelle's fragmentation at small times. In this direction,
we use the  following results due to Berestycki
\cite{Berestycki02}.

Let  $X_k(t)$  be the frequency of the $k$-th largest  block at time
$t$ of  Ruelle's  fragmentation. Recall that  Ruelle's
fragmentation can be constructed from a Poisson measure $K$ on
$[0,1[\times
\sr \times \mcn$ with intensity $(\nu_t(ds)dt) \otimes \sharp$.
Set
$$K=(t,S(t),k(t))_{t\in [0,1[}=(t,(s_1(t),s_2(t),\ldots),k(t))_{t\in [0,1[}.$$
Let $(S^{(i)}(t),t\in [0,1[
)=(s^{(i)}_1(t),s^{(i)}_2(t),\ldots,t\in [0,1[)$ be the Poisson measure
obtained from $K$  restricted to the atoms such that
$k(t)=i$. So, it is a Poisson measure with intensity $\nu_t(ds)dt$.

Set $$R(t)=\max_{s\le t} {s^{(1)}_2(s)}.$$

\begin{duge}\label{small}
\begin{itemize}
    \item For  $t$ small enough, we have
$X_1(t)=\exp(-\xi_t) a.s.$ where $\xi_t$ is an increasing process
with independent increments and such that :
$$\forall
k>0,\;\mce\left[\exp(-k
\xi_t)\right]=\frac{\Gamma\left(k+1-t\right)}{\Gamma\left(1-t\right)
\Gamma\left(k+1\right)}.$$
 \item $$X_2(t)\sim R(t) \mbox{ as } t\rightarrow 0^{+}
 \;a.s.$$
\end{itemize}
\end{duge}

Proof. The proof is the same as in Berestycki
\cite{Berestycki02}, since there, time-homogeneity of the
fragmentation plays no role. $\Box$

\vspace*{0.5cm}

Let us  now determine the behavior of $R(t)$.

\begin{prop}\label{encaRt}
Fix $T_0\in ]0,1/2[$. Then there exists three strictly positive
constants $C_1$,$C_2$, $C_3$  such that for all  $\lambda>0$ and
for all $t\in]0,T_0[$, $$ \exp(-C_1\lambda-C_3t)\le
\mcp\left(R(t)\le \exp\left(-\frac{\lambda}{t}\right)\right)\le
\exp(-C_2\lambda+C_3t).$$
\end{prop}

To estimate the distribution of  $R(t)$, we study $\nu_t(s_2\ge \varepsilon)$  for a fixed $\varepsilon$.
Indeed, $$\mcp(R(t)\le\varepsilon)=\exp(-\int_{0}^{t}\nu_u(s_2\ge\varepsilon)du),$$
and Proposition \ref{encaRt} follows from the following lemma :
\begin{duge}
Fix $T_0\in ]0,1/2[$. Then there exists three strictly positive
constants $C_1$,$C_2$, $C_3$  such that for all $\varepsilon \in
]0,1[$ and  for all $t\in]0,T_0[$,
$$-t(C_2\ln \varepsilon+C_3)\le \int_{0}^{t}\nu_u(s_2\ge\varepsilon)du\le -t(C_1\ln \varepsilon-C_3). $$
\end{duge}

Proof. We begin with the upper bound. If
$(s_i)_{i\ge 1}$ is an element of $\sr$, we denote
$(\tilde{s}_i)_{i\ge 1}$ a size-biaised rearrangement. We
have :
$$s_2\ge \varepsilon \Rightarrow s_1\le 1-\varepsilon \Rightarrow \tilde{s}_1\le 1-\varepsilon,$$

so
$$\nu_t(s_2\ge \varepsilon)\le \nu_t(s_1\le 1-\varepsilon) \le \nu_t( \tilde{s_1}\le 1-\varepsilon).$$

According to Theorem \ref{mesure}, we  know the law of
$\tilde{s}_1$ under $\nu_t$ :

\begin{eqnarray*}
\nu_t( \tilde{s}_1\le
1-\varepsilon)&=&\int_{0}^{1-\varepsilon}{(1-y)^{-1}y^{-t}dy}\\
&\le & \left(\int_{0}^{1/2}{2y^{-t}dy}
+\int_{1/2}^{1-\varepsilon}{2^{t}(1-y)^{-1}dy}\right)\\
&\le & \left(-2^{t}\ln \varepsilon
+\frac{2^t}{1-t} \right)\\
&\le & 2\left(-\ln \varepsilon
+2 \right) \hspace*{3cm}\mbox{ for } t\le \frac{1}{2}.\\
\end{eqnarray*}

So we obtain
$$ \int_{0}^{t}\nu_u(s_2\ge\varepsilon)du\le -t(2\ln
\varepsilon-4).$$

Let us now prove the lower bound. First, we will find a lower
bound for $\int_{0}^{t}\nu_u(\tilde{s}_2\ge\varepsilon)du$ and
then we will deduce the  lemma.

\begin{eqnarray*}
\nu_t(\tilde{s}_2\in dx)&=& \int_{0}^{1-x}\nu_t(\tilde{s}_1\in
dy)\nu_t(\tilde{s}_2\in dx\;|\; \tilde{s}_1\in dy)\\
&=& \frac{1
}{\Gamma(1-t)\Gamma(t)}\int_{0}^{1-x}
(1-y)^{-1}y^{-t}\left(\frac{x}{1-y}\right)^{-t}\left(1-\frac{x}{1-y}\right)^{t-1}\frac{dx}{1-y}dy\\
&=&\frac{ x^{-t}dx
}{\Gamma(1-t)\Gamma(t)}\int_{0}^{1-x}
(1-y)^{-1}y^{-t}\left(1-y-x\right)^{t-1}dy.\\
\end{eqnarray*}

Set
$$A= \int_{\varepsilon}^{1}
\int_{0}^{1-x}x^{-t} (1-y)^{-1}y^{-t}\left(1-y-x\right)^{t-1}dy
dx,$$ so $$\nu_t(\tilde{s}_2\ge\varepsilon)=\frac{1}
{\Gamma(1-t)\Gamma(t)}A.$$

We now calculate  a lower bound for $A$ :
\begin{eqnarray*}
A&=& \int_{0}^{1-\varepsilon}\int_{\varepsilon}^{1-y}x^{-t}
(1-y)^{-1}y^{-t}\left(1-y-x\right)^{t-1}dx dy\\
&=&
\int_{0}^{1-\varepsilon}\left(\int_{\frac{\varepsilon}{1-y}}^{1}z^{-t}
\left(1-z\right)^{t-1}dz\right) (1-y)^{-1}y^{-t}dy\\
&=&
\int_{\varepsilon}^{1}\left(\int_{\frac{\varepsilon}{y}}^{1}z^{-t}
\left(1-z\right)^{t-1}dz\right) y^{-1}(1-y)^{-t}dy \\
&\ge &
\int_{\varepsilon}^{1}\left(\int_{\frac{\varepsilon}{y}}^{1}
\left(1-z\right)^{t-1}dz\right) y^{-1}(1-y)^{-t}dy \\
&\ge &\frac{1}{t} \int_{\varepsilon}^{1}
\left(1-\frac{\varepsilon}{y}\right) y^{-1}(1-y)^{-t}dy \\
&\ge &\frac{1}{t} \int_{\varepsilon}^{1}
\left(1-\frac{\varepsilon}{y}\right) y^{-1}dy \\
&\ge& \frac{1}{t}\left(-\ln\varepsilon-1\right).
\end{eqnarray*}

So
$$\nu_t(\tilde{s}_2\ge\varepsilon)\ge\frac{1}{\Gamma(1-t)\Gamma(t)t}\left(-\ln\varepsilon-1\right).$$
As
 $\Gamma(1-t)\Gamma(t)t=\frac{\pi t}{\sin(\pi t)}$  is a positive function which is bounded on  $]0,T_0[$, let $1/C_2$ be its
 maximum. By integration, we obtain :

$$\int_{0}^{t}\nu_u(\tilde{s}_2\ge\varepsilon)du\ge t
C_2\left(-\ln\varepsilon-1\right).$$

We would like now to deduce the lower bound for
$\int_{0}^{t}\nu_u(s_2\ge\varepsilon)du$. We use
$$\nu_u(s_2\ge\varepsilon)\ge
\nu_u(\tilde{s}_2\ge\varepsilon)-\nu_u(\tilde{s}_2>
s_2),$$

and
$$\nu_u(\tilde{s}_2>
s_2)=\nu_u(\tilde{s}_2 =s_1)\le \nu_u(\tilde{s}_1 \neq
s_1)=\int_{\sr}(1-s_1)\nu_u(ds)\le
\int_{\sr}(1-\tilde{s_1})\nu_u(ds).$$
We have already seen that
$$\int_{\sr}(1-\tilde{s_1})\nu_u(ds)=\frac{1}{1-u} \hspace*{0.5cm}\mbox{ (cf. Formula (\ref{maj}))}.$$  So, for all
$t\le T_0$, we have
$$\int_{0}^{t}\nu_u(\tilde{s}_2>
s_2)du\le -\ln(1-t)\le \frac{1}{1-T_0}t.$$ Hence $$
\int_{0}^{t}\nu_u(s_2\ge\varepsilon)\ge t
\left(-C_2\ln\varepsilon-C_3\right).\Box$$

We can then deduce the lower-asymptotic behavior of  $X_2(t)$ from this theorem.

\begin{prop}\label{liminfrt}
There exists a constant
 $\delta>0 $ such that almost surely
$$\left\{
\begin{array}{lll}
\liminf_{t\rightarrow 0}|\ln t|^{\gamma/t}X_2(t)=0\;& \mbox{ if }\; \gamma<\delta\\
\liminf_{t\rightarrow 0}|\ln t|^{\gamma/t}X_2(t)=\infty \;&\mbox{ if }\;
 \gamma>\delta.
\end{array}
\right.$$
\end{prop}

Proof. According to Theorem \ref{small}, we just have to prove the
proposition replacing $X_2(t)$ by $R(t)$.  Set
$\gamma>\frac{1}{C_2}$. Choose $\beta>0$ such that
$\gamma>\frac{e^\beta}{C_2}$.
Set $t_i=e^{-i\beta}$ and $f(t)=\gamma \ln(-\ln t).$ For $t\in[0,e^{-1}[$, $f(t)$ is a decreasing positive function.  \\
For $t\in [t_{i+1},t_i]$, we have
$$R(t)\ge R(t_{i+1}) \mbox{ and } \exp\left(-\frac{f(t_i)}{t_i}\right) \ge \exp\left(-\frac{f(t)}{t}\right).$$

So if we prove
\begin{equation}\label{liminfr}
R(t_{i+1})\ge \exp\left(-\frac{f(t_i)}{t_i}\right)
\end{equation}
almost surely for $i$ large enough, then we will deduce

 $$ \forall \gamma>\frac{1}{C_2}\; \liminf_{t\rightarrow 0}(\ln\frac{1}{t})^{\gamma/t}R(t)\ge 1 \;\mbox{ a.s. }\; \mbox{ and so } \forall \gamma>\frac{1}{C_2}
 \;\liminf_{t\rightarrow 0}(\ln\frac{1}{t})^{\gamma/t}R(t)=\infty \;\mbox{ a.s. }\;$$

To prove (\ref{liminfr}), we apply Borel-Cantelli's Lemma.
 Using Proposition \ref{encaRt}, we obtain :

$$\mcp\left(R(t_{i+1})\le \exp\left(-\frac{f(t_i)}{t_i}\right)\right)\le K (\beta i)^{-C_2\gamma e^{-\beta} } .$$
Thanks to the choice of $\gamma$ and $\beta$, the serie converges.

For the second part of the proposition, we use an extension of
Borel-Cantelli's Lemma when the sum  diverges but the  events are
not
independent (cf. \cite{Kochen64}) :\\
Let $(H_i)_{i\ge 1}$ be a sequence of events such that $\sum
\mcp(H_i)$ diverges and
\begin{equation}\label{borel}
\forall N\ge 1, \;\frac{\sum_{i,j=1}^{N} \mcp(H_i\cap
H_j)}{\left(\sum_{i=1}^{N}\mcp(H_i)\right)^2}\le M.
\end{equation}
Then the set $\{i, \omega\in H_i\}$ is infinite with a  probability larger than  $1/M$.

In our case, we fix a  $\gamma<1/C_1$ and a  $\varepsilon>0$ such
that $(1+\varepsilon)\gamma C_1<1$. Set
$t_i=e^{-i^{1+\varepsilon}}$ and  $H_i=\{R(t_i)\le
(\ln(1/t_i))^{\gamma/t_i}\}$. Fix $i,j\ge 1$. Recall $R(t)$ is
the record process of a point Poisson process. So we have
\begin{eqnarray*}
\mcp(H_i\cap H_{j+i})&=&\mcp(H_i)\mcp(H_{i+j})\exp\left(\int_{0}^{t_{i+j}}\nu_u(s_2\ge (\ln(1/t_i))^{\gamma/t_i})du\right)\\
&\le& K \mcp(H_i)\mcp(H_{i+j})\exp\left((1+\varepsilon)C_1\gamma \ln i e^{-(1+\varepsilon)i^{\varepsilon}} \right)\\
&\le& K'\mcp(H_i)\mcp(H_{i+j}).
\end{eqnarray*}

(We have used  $(i+j)^{1+\varepsilon}-i^{1+\varepsilon}\ge
(1+\varepsilon)i^{\varepsilon}$ for all $i,j\ge 1$). With this
upper bound, we deduce that the sequence $H_i$ verifies
(\ref{borel}). We now have to prove that the sum of probabilities
diverges. Using  Proposition \ref{encaRt}, we obtain :

$$\sum_{i}\mcp(H_i)\ge K\sum_{i} i^{-C_1\gamma(1+\varepsilon)}.$$
Thus this series diverges thanks to our choice of $\gamma$ and
$\varepsilon$. We now apply the 0-1 law to prove that the
probability that  the set $\{i, \omega\in H_i\}$ is infinite equal
to 1.

So we have proved

$$\left\{
\begin{array}{lll}
\liminf_{t\rightarrow 0}(\ln\frac{1}{t})^{\gamma/t}R(t)=0\;& \mbox{ a.s. }\;& \forall \gamma<\frac{1}{C_1}\\
\liminf_{t\rightarrow 0}(\ln\frac{1}{t})^{\gamma/t}R(t)=\infty
\;&\mbox{ a.s. }\;& \forall \gamma>\frac{1}{C_2}.
\end{array}
\right.
$$

Thus we deduce that there exists almost surely a (random) critical
$\gamma_c\in]1/C_1,1/C_2[$ such that
$$\left\{
\begin{array}{lll}
\liminf_{t\rightarrow 0}(\ln\frac{1}{t})^{\gamma/t}R(t)=0\;&\;& \forall \gamma<\gamma_c\\
\liminf_{t\rightarrow 0}(\ln\frac{1}{t})^{\gamma/t}R(t)=\infty
\;&\;& \forall \gamma>\gamma_c.
\end{array}
\right.
$$

 By the 0-1 law, the law of $\gamma_c$ is trivial, i.e. it exists $\delta$ verifying
Proposition \ref{liminfrt}
 $\Box$

We can also determine the upper asymptotic behavior of $X_2(t)$ :

\begin{prop}\label{limsuprt}We have almost surely
$$\left\{
\begin{array}{lll}
\limsup_{t\rightarrow 0}\exp(\frac{1}{t}(-\ln(t))^{-\beta})X_2(t)=\infty\;& \mbox{ if } \; \beta>1\\
\limsup_{t\rightarrow 0}\exp(\frac{1}{t}(-\ln(t))^{-\beta})X_2(t)=0 \;&\mbox{ if } \; \beta\le1.\\
\end{array}
\right.
$$
\end{prop}

Proof. We use the same approach as for the infimum.
Fix $\beta>1$. Set $t_i=e^{-i}$ and  $f(t)=\exp(-\frac{1}{t}(-\ln(t))^{-\beta})$.
We want to prove that $R(t)\le f(t)$ almost surely for  $t$ small enough.
As $f$ is a decreasing function and R(t) an increasing process, we have $R(t)\le R(t_i)$ and $f(t_{i+1})\le f(t)$.
So we just have to prove that $R(t_i)\le f(t_{i+1})$ almost surely for $i$ large enough.

We have
\begin{eqnarray*}
\mcp\left(R(t_i)\ge f(t_{i+1})\right)&\le& 1-\exp\left(-C_3e^{-i}-C_1 e (i+1)^{-\beta}\right)\\
&\le& C_1 e i^{-\beta}+o(i^{-\beta}).
\end{eqnarray*}

This serie converges. So,  thanks to Borel-Cantelli's Lemma, we can conclude.

Let us now prove the case $\beta\le1$. Set  $t_i=e^{-i}$ and
$f(t)=\exp(-\frac{1}{t}(-\ln(t))^{-\beta})$. Set $H_i=\{R(t_i)\ge
f(t_i)\}$. Then we have
$$\sum_{i=1}^{N}\mcp(H_i)\ge  \sum_{i=1}^{N}\left(1-\exp\left(C_3e^{-i}-C_2 i^{-\beta}\right)\right).$$

The right term is equivalent to $\sum_{i=1}^{N}C_2i^{-\beta}$, so  it diverges.

We have now to check  the condition (\ref{borel}) to apply the generalized
Borel-Cantelli's Lemma.

\begin{eqnarray*}
\mcp(H_i\cap H_{i+j})&=& 1-\mcp(\overline{H_i})-\mcp(\overline{H_{i+j}})+\mcp(\overline{H_i}\cap\overline{H_{i+j}} )\\
 &=&  1-\mcp(\overline{H_i})-\mcp(\overline{H_{i+j}})+\mcp(\overline{H_i})\mcp(\overline{H_{i+j}})\exp\left(\int_{0}^{t_{i+j}}\nu_u(s_2\ge f(t_i))du\right)\\
&\le&
\mcp(H_i)\mcp(H_{i+j})+\exp\left(\int_{0}^{t_{i+j}}\nu_u(s_2\ge
f(t_i))du\right)-1.
\end{eqnarray*}

Then remark

$$\exp\left(\int_{0}^{t_{i+j}}\nu_u(s_2\ge f(t_i))du\right)\le \exp\left(C_3 e^{-i-j}+C_1i^{-\beta}e^ {-j}\right).$$

Hence we deduce

$$\sum_{i,j=1}^{N}\left( \exp\left(\int_{0}^{t_{i+j}}\nu_u(s_2\ge f(t_i))du\right)-1 \right)\le K\sum_{i=1}^{N}i^{-\beta}.$$

So

$$\frac{\sum_{i,j=1}^{N}\left( \exp\left(\int_{0}^{t_{i+j}}\nu_u(s_2\ge f(t_i))du\right)-1 \right)}{\sum_{i=1}^{N}{\mcp(H_i)}}$$ is bounded and thus
the condition (\ref{borel}) is true.

So, we can conclude for the case $\beta<1$.
For $\beta=1$, we just have
$$\limsup_{t\rightarrow 0}R(t)\exp\left(-\frac{1}{t\ln t}\right)\le1 \;\mbox{ a.s. }$$
Remark then that the same demonstration works with $\gamma f(t)$ instead of $f(t)$ with  $\gamma$ positive constant.
So, we have
$$\limsup_{t\rightarrow 0}R(t)\exp\left(-\frac{1}{t\ln t}\right)\le\gamma \;\mbox{ a.s. }$$
and thus

$$\limsup_{t\rightarrow 0}R(t)\exp\left(-\frac{1}{t\ln t}\right)=0 \;\mbox{ a.s. } \Box$$

 \vspace*{2cm}
~\nocite{Bertoin00} ~\nocite{Bertoinlegall00} ~\nocite{Bertoin01}
~\nocite{Marchal1} ~\nocite{Marchal2} ~\nocite{Ruelle87}
~\nocite{BovierKurkova041} ~\nocite{BovierKurkova044}
~\nocite{Durrett91}
~\nocite{BovierKurkova03}

\bibliographystyle{plain}
\bibliography{Ruellescascades}

\end{document}